\documentclass[final,leqno,onefignum,onetabnum]{siamltex1213}

\usepackage{hyperref}
\usepackage{amsmath}
\usepackage{amssymb} 
\usepackage{algorithm}
\usepackage{algorithmic}
\usepackage{multirow}
\usepackage{pstool}
\usepackage{multicol}
\allowdisplaybreaks
\allowdisplaybreaks

\newcommand{\func}[1]{#1}
\newcommand{\oper}[1]{#1}

\newcommand{\trace}{\operatorname{tr}}
\newcommand{\Var}{\operatorname{Var}}

\newcommand{\bE}{{\mathbb{E}}}

\newcommand{\bR}{{\mathbb{R}}}

\newcommand{\bS}{{\mathbb{S}}}

\newcommand{\cB}{\mathcal{B}}
\newcommand{\cC}{\mathcal{C}}

\newcommand{\cH}{\mathcal{H}}
\newcommand{\cJ}{\mathcal{J}}

\newcommand{\cL}{\mathcal{L}}
\newcommand{\cM}{\mathcal{M}}
\newcommand{\cN}{\mathcal{N}}
\newcommand{\cP}{\mathcal{P}}
\newcommand{\cR}{\mathcal{R}}

\newcommand{\cU}{\mathcal{U}}
\newcommand{\cV}{\mathcal{V}}

\newcommand{\cX}{\mathcal{X}}
\newcommand{\cY}{\mathcal{Y}}

\newcommand{\cZ}{\mathcal{Z}}

\newcommand{\bstau}{\boldsymbol{\tau}}

\newcommand{\bse}{{\boldsymbol{e}}}

\newcommand{\bsn}{{\boldsymbol{n}}}

\newcommand{\bsv}{\boldsymbol{v}}
\newcommand{\bsu}{{\boldsymbol{u}}}

\newcommand{\bssigma}{{\boldsymbol{\sigma}}}

\newcommand{\beq}{\begin{equation}}
\newcommand{\eeq}{\end{equation}}
\newcommand{\ba}{\begin{array}}
\newcommand{\ea}{\end{array}}

\usepackage{framed,multirow}

\usepackage{latexsym}

\usepackage{url}
\usepackage{xcolor}
\definecolor{newcolor}{rgb}{.8,.349,.1}
\title{Taylor approximation and variance reduction for PDE-constrained optimal control under uncertainty 
} 

\author{Peng Chen \thanks{Institute for Computational Engineering \& Sciences, The University of Texas at Austin, Austin, TX 78712 (\email{peng@ices.utexas.edu, uvilla@ices.utexas.edu}).}  \and Umberto Villa \footnotemark[2] \and Omar Ghattas \thanks{Institute for Computational Engineering \& Sciences, Department of Mechanical Engineering, and Department of Geological Sciences, The University of Texas at Austin, Austin, TX 78712 (\email{omar@ices.utexas.edu}).} }
\begin{document}
\maketitle
\slugger{sisc}{xxxx}{xx}{x}{x--x}

\begin{abstract}
In this work we develop a scalable computational framework for the solution of PDE-constrained optimal control under high-dimensional uncertainty. Specifically, we consider a mean-variance formulation of the control objective and employ
a Taylor expansion with respect to the uncertain parameter either to directly approximate the control objective or as a control variate for variance reduction. The expressions for the mean and variance of the Taylor approximation are known analytically, although their evaluation requires efficient computation of the trace of the (preconditioned) Hessian of the control objective. 
We propose to estimate this trace by solving a generalized eigenvalue problem using a randomized algorithm that only requires the action of the Hessian on a small number of random directions. Then, the computational work does not depend on the nominal dimension of the uncertain parameter, but depends only on the  effective dimension (i.e., the rank of the preconditioned Hessian), thus ensuring scalability to high-dimensional problems. 
Moreover, to increase the estimation accuracy of the mean and variance of the control objective by the Taylor approximation, we use it as a control variate for variance reduction, which results in considerable computational savings (several orders of magnitude) compared to a plain Monte Carlo method.
In summary, our approach amounts to solving an optimal control constrained by the original PDE and a set of linearized PDEs, which arise from the computation of the gradient and Hessian of the control objective with respect to the uncertain parameter. 
We use the Lagrangian formalism to derive expressions for the gradient with respect to the control and apply a gradient-based optimization method to solve the problem. 
We demonstrate the accuracy, efficiency, and scalability of the proposed computational method for two examples with high-dimensional uncertain parameters: subsurface flow in a porous medium modeled as an elliptic PDE, and turbulent jet flow modeled by the Reynolds-averaged Navier--Stokes equations coupled with a nonlinear advection-diffusion equation characterizing model uncertainty. In particular, for the latter more challenging example we show scalability of our algorithm up to one million parameters resulting from discretization of the uncertain parameter field.

\end{abstract}



\pagestyle{myheadings}
\thispagestyle{plain}
\markboth{Taylor approximation and variance reduction for optimal control under uncertainty}{P. Chen, U. Villa, and O. Ghattas}

\section{Introduction}

Optimal control (and more general optimization) constrained by partial differential equations (PDEs) are ubiquitous in many applications of practical relevance for science and engineering. The objective of such problems is to seek an optimal control that minimizes a cost functional, which often consists of a control objective and a penalty term. The control objective is often given as a distance between the solution of the PDE and a desired state or more generally formulated as a quantity of interest related to the solution, while the penalty term reflects the cost of the control or imposes regularity of the control. In practical applications, uncertainties in the PDE model are inevitable and can arise from various sources, such as the PDE coefficients, initial or boundary conditions, external sources, or geometries. Different realizations of these uncertain parameters may lead to significant differences in the optimal control. Mathematical theories and computational methods have been developed for several decades to deal with deterministic PDE-constrained optimal control 
\cite{Lions1971, GlowinskiLions1996, Gunzburger2003, BieglerGhattasHeinkenschlossEtAl03, HinzePinnauUlbrichEtAl2008, Troeltzsch2010, BorziSchulz2011}. 
More recently, PDE-constrained optimal control under uncertainty has become an important area of research and has drawn increasing attention \cite{BorziSchulzSchillingsEtAl2010, SchillingsSchmidtSchulz2011, HouLeeManouzi2011,GunzburgerLeeLee2011,RosseelWells2012,
KouriHeinkenschloosVanBloemenWaanders2012,TieslerKirbyXiuEtAl2012,ChenQuarteroniRozza2013,LassilaManzoniQuarteroniEtAl2013, ChenQuarteroni2014, KouriHeinkenschlossRidzalEtAl14, KunothSchwab2013, NgWillcox2014, ChenQuarteroniRozza2016, KunothSchwab2016, KouriSurowiec2016, BennerOnwuntaStoll2016, AlexanderianPetraStadlerEtAl17, AliUllmannHinze2017}.
  
Different formulations of the cost functional have been studied to incorporate uncertainty into the PDE-constrained optimal control. A straightforward choice is to optimize the mean of the control objective, i.e., the integration of the control objective with respect to the probability measure of the uncertain parameter \cite{BorziSchulzSchillingsEtAl2010,HouLeeManouzi2011,GunzburgerLeeLee2011,
KouriHeinkenschloosVanBloemenWaanders2012, ChenQuarteroni2014}. However, this formulation does not control the variability
of the control objective. To avoid a possible and undesired large variation of the control objective, which may represent the risk of a system failure, one can include the variance or higher moments of the control objective in the cost functional \cite{RosseelWells2012, BennerOnwuntaStoll2016, TieslerKirbyXiuEtAl2012}. An alternative approach is use the Value-at-Risk (VaR) or conditional VaR \cite{KouriSurowiec2016}, which measure the (conditional) probability of the control objective surpassing a certain critical value. An extreme choice is robust optimization (also referred as min-max optimization), where one seeks a control that yields the optimal value of the objective in the worst-case scenario, i.e., the extreme value of the control objective in the range of the uncertain parameter \cite{bertsimas2011theory, LassilaManzoniQuarteroniEtAl2013}.

A common challenge in all of these formulations is to compute the statistics of the control objective (i.e., mean, variance, higher moments, or conditional probabilities); this involves the integration of the control objective or related quantities with respect to the probability measure of the uncertain parameter. This challenge becomes more severe for high-dimensional uncertain parameter space, which is often the case when the parameter is a uncertain spatially correlated field. A straightforward approach is to apply a Monte Carlo method to compute the statistics by taking the average of the control objective or related quantities at a set of samples randomly drawn according to the probability measure of the parameter. However, the converge rate of the Monte Carlo estimator is just $O(M^{-1/2})$, where $M$ is the number of samples. This may lead to an extremely large number of samples to achieve a target accuracy, and therefore may not be feasible when the state equation is a complex nonlinear PDE. Quasi Monte Carlo methods using low-discrepancy sequences improve the convergence rate to $O((\log(M))^d/M)$ for a $d$-dimensional parameter \cite{DickKuoSloan2013}, and is thus attractive for low to medium dimensional problems.
Stochastic Galerkin and stochastic collocation methods have also been used to compute the statistical moments in the cost functional \cite{HouLeeManouzi2011, GunzburgerLeeLee2011, 
TieslerKirbyXiuEtAl2012, RosseelWells2012, KouriHeinkenschloosVanBloemenWaanders2012, ChenQuarteroniRozza2013, LassilaManzoniQuarteroniEtAl2013,
ChenQuarteroniRozza2013, ChenQuarteroni2014,
KunothSchwab2013, KunothSchwab2016}, provided that a suitable finite dimensional parametrization of the uncertain parameter, such as a truncated Karhunen--Lo\`eve expansion, is applicable. These methods achieve fast convergence when the control objective depends smoothly on the low-dimensional parameter, but suffer from the so-called \emph{curse of dimensionality}, i.e., the convergence rate quickly deteriorates as the dimension of the parameter increases. Recent advances in adaptive and anisotropic sparse quadrature \cite{SchillingsSchwab2013, Chen2018} and high-order quasi Monte Carlo methods \cite{DickLeGiaSchwab2016} have been shown to achieve a dimension-independent convergence rate $N^{-s}$ with $s$ potentially much larger than $1/2$. Here $s$ is independent of the nominal dimension of the uncertain parameter and depends only on the choice of parametrization and on the regularity of the control objective with respect to the parameter, thus mitigating or even breaking the curse of dimensionality. However, the convergence may still be slow if the choice of the parametrization does not correctly capture the effective parameter dimensions.

Another computational challenge is that, in applications of practical relevance, the PDE may lead, after discretization, to possibly nonlinear large-scale systems that are extremely expensive to solve. Thus, only a limited number of high-fidelity PDE solves can be afforded. This challenge prevents a direct application of most of the computational methods introduced above as they require a large number of evaluations of the control objective for computing its statistics. To tackle this challenge, reduced basis methods \cite{RozzaHuynhPatera2008, QuarteroniManzoniNegri2015, HesthavenRozzaStamm2015} or other model reduction techniques \cite{BennerGugercinWillcox2015} can be employed to exploit the low dimensionality of the solution manifold, even when the parameter lives in a high-dimensional space. Such methods solve the high-fidelity state equation at a few carefully selected samples of the parameter, and use such high-fidelity solutions (known as snapshots) as the basis to approximate the solution at any other realization of the uncertain parameter by a (Petrov)-Galerkin projection. Solving the projected problem in the reduced basis space allows for considerable computational saving compared to solving in the full space, such as a finite element space with a large number of degrees of freedom. These techniques have been successfully developed for solving parametric and stochastic PDE-constrained optimal control with several classes of PDE models \cite{NegriRozzaManzoniEtAl2013, KaercherGrepl2014, BennerSachsVolkwein2014, NegriManzoniRozza2015, ChenQuarteroni2014, ZahrFarhat2015, ChenQuarteroniRozza2016, BaderKaercherGreplEtAl2016, ChenQuarteroniRozza2017}. Nevertheless, difficulties in using such techniques arise for highly nonlinear problems that require effective affine approximation (see \cite{QuarteroniManzoniNegri2015, HesthavenRozzaStamm2015} and references therein) or when the solution manifold becomes high-dimensional, even if the effective dimension of the manifold of the control objective is low \cite{Bui-ThanhWillcoxGhattasEtAl2007, ChenGhattas2017}. 

In this work, we propose a scalable computational framework that uses Taylor expansions of the control objective with respect to the uncertain parameter to efficiently and accurately approximate the cost functional. 
This framework takes advantage both of the smooth dependence of the control objective on the uncertain parameter and of the low-rank structure of the Hessian of the control objective, which exposes the low effective dimension of the nominally high-dimensional parameter space.
Such Taylor expansions can be used directly to approximate the moments of the control objective as proposed in \cite{AlexanderianPetraStadlerEtAl17} or as control variates to reduce the variance of Monte Carlo estimators.
Analytic expressions for the mean and variance of the Taylor expansions are known; we refer to \cite{AlexanderianPetraStadlerEtAl17} for such expressions in the case of linear and quadratic approximations. These expressions depend on the trace of the preconditioned Hessian (second order derivative of the control objective with respect to the uncertain parameter) and/or its square. While for low-dimensional problems it is possible to explicitly compute the Hessian matrix (see \cite{GhateGiles2007, GhattasBark97}), this approach is not tractable for high-dimensional problems, where the Hessian is a formally large dense operator implicitly defined by its action on a given vector.  Gaussian trace estimators (see \cite{AvronToledo2011}) only need to perform the action of the Hessian on a number of random directions, and have been successfully used for the solution of optimal experimental design and optimization under uncertainty in \cite{AlexanderianPetraStadlerEtAl14}. 
However, such estimators may still require a large number of Hessian actions to achieve high accuracy, see \cite{AvronToledo2011} for probablistic bounds on the number of Hessian actions necessary to achieve a target accuracy with high probability.
In this work, we approximate the trace of the preconditioned Hessian with the sum of its dominant eigenvalues. Randomized algorithms \cite{HalkoMartinssonTropp2011,SaibabaLeeKitanidis2016} are computationally efficient methods to compute dominant eigenvalues of a large-scale implicitly defined operator whose spectra exhibit fast decay. The spectral properties of the Hessian operator have been extensively studied, and it has been observed numerically or proven analytically that for many problems the Hessian operator is either nearly low-rank or its eigenvalues exhibit fast decay \cite{BashirWillcoxGhattasEtAl08, FlathWilcoxAkcelikEtAl11,
  Bui-ThanhGhattas12a, Bui-ThanhGhattas13a, Bui-ThanhGhattas12,
  Bui-ThanhBursteddeGhattasEtAl12_gbfinalist,
  Bui-ThanhGhattasMartinEtAl13, ChenVillaGhattas2017,
  AlexanderianPetraStadlerEtAl16, AlexanderianPetraStadlerEtAl17,
  AlexanderianPetraStadlerEtAl14, CrestelAlexanderianStadlerEtAl17,
  PetraMartinStadlerEtAl14, IsaacPetraStadlerEtAl15,
  MartinWilcoxBursteddeEtAl12, Bui-ThanhGhattas15, ChenGhattas18b}.
Thus, this approach can provide highly accurate approximation of the trace at a cost (measured in number of Hessian actions) much smaller than Gaussian trace estimators. 

The Taylor approximation may be not sufficiently accurate and introduces bias in the estimation of the mean and variance of the control objective. 
Specifically, the norm of the remainder of the Taylor expansion scales as $O(\text{tr}(\cC)^{(r+1)/2})$ where $\cC$ is the covariance of the uncertain parameter and $r$ is the order of the truncated Taylor expansion, e.g., $r= 1$ and $2$ for linear and quadratic approximations, respectively \cite{AlexanderianPetraStadlerEtAl17}. 
To address this issue, we propose to use the Taylor approximation as a control variate to reduce the variance of the Monte Carlo estimator for the mean and variance of the control objective, or use Monte Carlo correction for the remainder of the Taylor expansion from another viewpoint. 
In cases when the Taylor approximation is highly correlated with the exact control objective, reduction in the computational cost can be dramatic, up to several orders of magnitude in our numerical examples. 

We consider two examples of PDE-constrained optimal control under uncertainty to numerically illustrate the scalability, accuracy and efficiency of our computational framework. The first is a problem of optimal control of subsurface flow in a porous medium, where the control is the injection rate at given well locations, the state equation is an elliptic PDE representing single phase flow with uncertain permeability, and the objective is to drive the pressure field to a desired state. The second is an optimal boundary control for a turbulent jet flow, where the control is the inlet velocity profile, the state equations are the Reynolds-averaged Navier--Stokes equations coupled with a nonlinear stochastic advection-diffusion equation representing model uncertainty, and the objective is to maximize the jet width at a cross-section a certain distance from the inlet boundary. In both cases, the uncertain parameter is infinite-dimensional in the continuous setting and becomes high-dimensional after discretization. 

For solution of the PDE-constrained optimal control with Taylor approximation and variance reduction we employ a gradient-based method, specifically the limited memory BFGS method \cite{NocedalWright2006}. The cost functional involves the analytical expressions for the mean and variance of the Taylor approximation and the Monte Carlo correction when using the Taylor approximation for variance reduction; the constraint involves the state (possibly nonlinear) PDE and a set of linearized PDEs used in the computation of the gradient and Hessian of the control objective with respect to the uncertain parameter.
We use the Lagrangian formalism to derive the gradient of the cost functional with respect to the control. Specifically, we compute mixed high-order derivatives of the Lagrangian functional with respect to the state, the adjoint, the uncertain parameter, and the control by symbolic differentiation in a variational setting before discretization. This allows us to efficiently apply our framework to rather complex PDE models such as the turbulent jet example, which features coupled nonlinear PDEs, nonlinear numerical stabilization, and weak imposition of Dirichlet boundary conditions.
We demonstrate numerically that our method scales independent of the uncertain parameter dimension, in that the number of required PDE solves does not depend on the discretization of the parameter field, and instead depends only on the effective or instrinsic dimension of the problem, i.e., the number of parameter dimensions to which the control objective is sensitive.

The rest of the paper is organized as follows: We present the basic setting for the PDE-constrained optimal control under uncertainty in Section \ref{sec:setting} using a general notation for the uncertain parameter, the PDE model, the cost functional, and the optimal control problem. We also present a sample average approximation of the mean and variance of the control objective in this section. Section \ref{sec:taylor} is devoted to the development of the Taylor approximation and variance reduction for the evaluation of the mean and variance of the control objective. We also present computation of the gradient and Hessian of the control objective with respect to the uncertain parameter and a double-pass randomized algorithm for solution of generalized eigenvalue problems. This is followed by Section \ref{sec:gradient} where we derive the gradient of the cost functional with respect to the control for different approximations including the sample average approximation, Taylor approximation, and Taylor approximation with Monte Carlo correction. Numerical results for two examples are presented in Section \ref{sec:numerics} for the demonstration of the scalability, accuracy, and efficiency of the proposed computational methods. Section \ref{sec:conclusion} closes with some conclusions and perspectives.

\section{PDE-constrained optimal control under uncertainty}
\label{sec:setting}

In this section we present the PDE-constrained optimal control under uncertainty in an abstract setting.
To this aim, let us first introduce the notation used throughout the paper. Let $\cX$ be a separable Banach space and $\cX'$ the dual space, then $_\cX \langle \cdot, \cdot \rangle_{\cX'}$ denotes the duality pairing between the spaces $\cX$ and $\cX'$. For ease of notation, we will omit to specify the subscritps $\cX$ and $\cX'$ and simply write $\langle \cdot, \cdot \rangle$ when the spaces can be inferred from the context without ambiguity. Given two separable Banach spaces $\cX$ and $\cY$ and a map $f: \cX \times \cY \mapsto \mathbb{R}$, $\partial_x f(x, y) \in \cX'$ denotes the Fr\'echet derivative of $f(x,y)$ with respect to $x$ evaluated at $(x,y)$, which satisfies
\beq
\lim_{\tilde{x} \to 0}|f(x+\tilde{x}, y) - f(x, y) - \langle \tilde{x},  \partial_x f(x, y) \rangle|/||\tilde{x}||_\cX = 0.
\eeq
Let $\partial_{xy}f(x,y): \cY \mapsto \cX'$ denote the Fr\'echet derivative of $\partial_x f(x,y)$ with respect to $y$ evaluated at $(x,y)$, or the second order (mixed) Fr\'echet derivative of $f(x,y)$ with respect to $x$ and $y$ evaluated at $(x,y)$. Similarly, $\partial_{yx} f(x,y): \cX \mapsto \cY'$ denotes the Fr\'echet derivative of $\partial_y f(x,y)$ with respect to $x$ evaluated at $(x,y)$, which is the adjoint operator of $\partial_{xy}f(x,y)$ and satisfies 
\beq\label{eq:adjoint_op}
\langle \tilde{x}, \partial_{xy} f \, \hat{y} \rangle: = \;_\cX \langle \tilde{x}, \partial_{xy} f \, \hat{y} \rangle_{\cX'} = \;_\cY \langle \hat{y}, \partial_{yx} f \, \tilde{x} \rangle_{\cY'} = :\langle \hat{y}, \partial_{yx} f \, \tilde{x} \rangle, \; \forall \tilde{x} \in \cX, \, \hat{y} \in \cY.
\eeq
where we have omitted the argument $(x,y)$ for simplicity.
Finally, let us introduce the separable Banach spaces $\cU, \cV, \cM,$ and $\cZ$ in which the state $u$, the adjoint $v$, the uncertain parameter $m$, and the deterministic control $z$ live, respectively. The exact definition of these spaces is problem specific and is specified in Sec.\ \ref{sec:numerics}.


\subsection{The uncertain parameter}
We assume that the uncertain parameter $m \in \cM$ has mean $\bar{m} \in \cM$ and covariance $\cC:\cM \to \cM'$. The computational framework developed in this work is applicable as long as this assumption is satisfied regarding the uncertain parameter, which can be both finite and infinite dimensional. In particular, to demonstrate the scalability of our computational methods, we consider an infinite-dimensional Gaussian random field with distribution $\mu = \cN(\bar{m}, \cC)$ for the uncertain parameter, which is one of the most important building blocks in modern spatial statistical modeling. A general covariance structure is known as Mat\'ern covariance. It is equivalent to that the uncertain parameter $m$ is given as the solution of some linear fractional stochastic PDE \cite{LindgrenRueLindstroem11}. In this work we consider this type of covariance, i.e., 
\beq\label{eq:covariance}
\cC = (-\alpha_1 \nabla \cdot (\Theta \nabla) +\alpha_2 I)^{-\alpha},
\eeq
where $\nabla \cdot$, $\nabla$, and $I$ are the divergence, gradient, and identity operators, respectively; the parameter $\alpha > 0$ reflects the smoothness of the random field, which renders the covariance of trace class when $\alpha > d/2$ for spatial dimension $d = 1, 2, 3$, $\alpha_1, \alpha_2 > 0$ controls its correlation length and variance, while $\Theta \in \bR^{d\times d}$ enables its spatial anisotropy.

\subsection{The PDE model and control objective}

We consider a state problem modeled by the following PDE model. Given a realization of the uncertain parameter $m \in \cM$ and a deterministic control $z \in \cZ$, find the state $u\in \cU$, such that
\beq\label{eq:PDE_strong}
\cR(u, m, z) = 0 \quad \text{ in }  \cV',
\eeq
where $\cR(\cdot, m, z) : \cU  \to \cV'$ denotes a (possibly nonlinear) operator from $\cU$ to the dual of $\cV$.
We then define the weak form associated to the state equation \eqref{eq:PDE_strong} by means of duality pairing, and write
\beq\label{eq:PDE} 
r(u, v, m, z) := \; _{\cV} \langle v, \cR(u, m, z) \rangle_{\cV'} = 0 \quad \forall v \in \cV.
\eeq 
Note that $r(u, v, m, z)$ is linear with respect to the adjoint $v$ and possibly nonlinear with respect to the state $u$. 

As for the control objective, we consider a real-valued, possibly nonlinear, functional $Q$ of the state $u\in \cU$,
\beq
Q: \cU \to \bR.
\eeq 
For simplicity, we assume that $Q$ depends on $m$ and $z$ only implicitly through the solution of the state equation \eqref{eq:PDE}.
The general case where $Q$ has an explicit dependence on $m$ and $z$ can be treated with slight modifications in the derivation of the cost functional and its gradient, as we address along the presentation in Sec.\ \ref{sec:taylor} and Sec.\ \ref{sec:gradientopt}.

\subsection{The PDE-constrained optimal control under uncertainty} 

Since the control objective $Q$ (implicitly) depends on the uncertain parameter $m$ through the state $u$, we face PDE-constrained optimal control under uncertainty.
In this work we consider a mean-variance cost functional of the form
\beq\label{eq:objective}
\cJ(z) = \bE[Q] + \beta \Var[Q] + \cP(z),
\eeq 
which is often used in stochastic programming \cite{ShapiroDentchevaRuszczynski09}.
Above $\bE[Q]$ and $\Var[Q]$ are, respectively, the mean and variance of $Q$ with respect to the probability measure $\mu$ of the uncertain parameter $m$, and are defined as
\beq
\bE[Q] := \int_\cM Q(u(m,z)) d\mu(m),
\eeq
and 
\beq
\Var[Q] := \int_\cM (Q(u(m,z)) - \bE[Q])^2 d\mu(m) = \bE[Q^2] - (\bE[Q])^2.
\eeq
As the control objective $Q$ (implicitly) depends on the control $z$ through the state $u$, so do $\bE[Q]$ and $\Var[Q]$.
The scaling factor $\beta \geq 0$ in \eqref{eq:objective} controls the relative weight of the mean and variance, and $\cP(z)$ is a penalty term. 

Recalling that the state $u$ (implicitly) depends on the uncertain parameter $m$ and control $z$ through the solution of the PDE \eqref{eq:PDE},
we then formulate the PDE-constrainted optimal control under uncertainty as: find $z^* \in \cZ$ such that  
\beq\label{eq:optimization}
z^* = \operatorname{arg}\,\min_{z \in \cZ} \cJ(z), \text{ subject to the problem } \eqref{eq:PDE}.
\eeq

To evaluate the mean and variance in the cost functional \eqref{eq:objective}, a plain Monte Carlo sample average approximation can be used. We briefly discuss this approach below, while in the following section we introduce more sophisticated approximation that allow for more efficient algorithms by exploiting high-order derivatives.

\subsection{Sample average approximation}
\label{sec:MCinteg}
The (plain Monte Carlo) sample average approximation of the mean reads 
\beq
\bE[Q] \approx  \widehat{Q} := \frac{1}{M_1} \sum_{i=1}^{M_1} Q(m_i),
\eeq
where $m_i$, $i = 1, \dots, M_1$, are independent indentically distributed (i.i.d.) random samples drawn from the probability measure $\mu$.
Similarly, the sample average approximation of the variance is given by
\beq
\Var[Q] \approx \widehat{V}_Q := \frac{1}{M_2} \sum_{i=1}^{M_2} Q^2(m_i) - \left(\frac{1}{M_2} \sum_{i=1}^{M_2} Q(m_i)\right)^2,
\eeq
where $m_i$, $i = 1, \dots, M_2$, are i.i.d. random samples from $\mu$.
Since each evaluation of $Q(m_i)$ involves solving the possibly nonlinear state equation \eqref{eq:PDE}, it is common practice to choose $M_1 = M_2 = M$ and reuse the same samples $m_i$ to estimate both the mean and variance of $Q$.
However, a more careful choice of $M_1$ and $M_2$ relies on estimating the variance of $Q$ and $Q^2$ and aims at minimizing the mean square error of $\widehat{Q} + \beta\widehat{V}_Q$.
{Note that a large number of samples is required to obtain a relatively accurate approximation of $\bE[Q]$ and $\Var[Q]$ due to the slow convergence of Monte Carlo sampling, thus rendering this approach computationally prohibitive for practical applications due to the large number of (nonlinear) PDE solves at each optimization step.}

\section{Taylor approximation and variance reduction}
\label{sec:taylor}
We develop the main approximation methods in this section, which include the approximation by (low-order) Taylor expansion of the control objective, the computation of the gradient and Hessian of the control objective with respect to the uncertain parameter used in the Taylor expansion, the randomized algorithm to solve a generalized eigenvalue problem for trace estimation, and the Monte Carlo correction for the remainder of the Taylor expansion.

\subsection{Taylor approximation of the control objective}
In this section we derive expressions for the Taylor expansion of the control objective $Q(u(m,z))$ with respect to the uncertain parameter $m$ for a fixed value of the control $z$. In what follows, we denote with the subscript $m$ the total derivative of $Q(u(m,z))$ with respect to $m$ and we assume sufficient smoothness so that all Fr\'echet derivatives are well-defined. For ease of notation, we also denote with $\langle \cdot, \cdot \rangle$ the duality-pairing $_\cM \langle \cdot, \cdot \rangle_{\cM'}$, and, with a slight abuse of notation, we write $Q(m)$ to indicate the control objective $Q(u(m,z))$ evaluated at a fixed value of $z$.

A formal Taylor expansion of the control objective evaluated at $\bar{m} \in \cM$ is given by 
\beq
Q(m) = Q(\bar{m}) + \langle m - \bar{m}, \func{Q}_m(\bar{m}) \rangle + \frac{1}{2} \langle m-\bar{m}, \oper{Q}_{mm}(\bar{m})\,(m-\bar{m})\rangle + \cdots,
\eeq
where $\func{Q}_m(m_0) \in \cM'$ and $\oper{Q}_{mm}(\bar{m}): \cM \to \cM'$ are the first and second order Fr\'echet derivatives of $Q$ 
with respect to $m$ evaluated at the nominal $\bar{m} \in \cM$. In what follows, we refer to such derivatives as the gradient and Hessian of $Q$, respectively. Based on the Taylor expansion, we can define the linear approximation of $Q$ as 
\beq\label{eq:linearQ}
Q_{\text{lin}}(m) = Q(\bar{m}) + \langle m - \bar{m}, \func{Q}_m(\bar{m})\rangle,
\eeq
and the quadratic approximation of $Q$ as
\beq\label{eq:quadraticQ}
Q_{\text{quad}}(m) = Q(\bar{m}) + \langle m - \bar{m}, \func{Q}_m(\bar{m})\rangle + \frac{1}{2} \langle m-\bar{m}, \oper{Q}_{mm}(\bar{m})\,(m-\bar{m}) \rangle.
\eeq
We remark that $Q_{\text{lin}}(m)$ and $Q_{\text{quad}}(m)$ are polynomials with respect to the uncertain parameter $m$, and depend on the control $z$ through the solution of the state equation \eqref{eq:PDE} at $m = \bar{m}$.

For the linear approximation, given the mean $\bar{m}$ and covariance $\cC$ of the uncertain parameter $m$, we have
\beq\label{eq:linearE}
\bE[Q_{\text{lin}}] = Q(\bar{m}), \text{ and } \Var[Q_{\text{lin}}] = \langle  \func{Q}_m(\bar{m}), \cC \func{Q}_m(\bar{m}) \rangle.
\eeq
In fact, $Q_{\text{lin}}(m)$ is Gaussian if the uncertain parameter $m \sim \cN(\bar{m}, \cC)$.
For the quadratic approximation, we have
\beq\label{eq:quadraticE}
\bE[Q_{\text{quad}}] = \bE[Q_{\text{lin}}] + \frac{1}{2} \trace(\cH), \text{ and } \Var[Q_{\text{quad}}] = \Var[Q_{\text{lin}}]+ \frac{1}{2} \trace(\cH^2),
\eeq
where $\cH = \cC^{1/2} \oper{Q}_{mm}(\bar{m}) \cC^{1/2}$ is the covariance preconditioned Hessian, and $\trace(\cdot)$ denotes the trace operator, see details in \cite{AlexanderianPetraStadlerEtAl17}. Efficient algorithms to accurately estimate the trace of $\cH$ and $\cH^2$ will be discussed in Sec.\ \ref{subsec:randtrace}.

Finally, following \cite{AlexanderianPetraStadlerEtAl17}, we approximate the cost functional $J(z)$ in \eqref{eq:objective} by replacing $\bE[Q]$ and $\text{Var}[Q]$ with the moments of the linear or quadratic Taylor expansions of the control objective.

\subsection{Computation of the gradient and Hessian of the control objective with respect to the uncertain parameter}
\label{subsec:gradientHessian}
In this section we use the Lagrangian formalism to derive the expressions for the gradient $Q_m$ and for the action of the Hessian $Q_{mm}$ on a given direction $\hat{m}$ evaluated at $\bar{m}$. 
In what follows, for ease of notation, we use the overline to indicate quantities evaluated at $m = \bar{m}$ and define
\beq\label{eq:r0Q0}
\bar{r} := r(u,v,\bar{m},z), \text{ and } \bar{Q} := Q(u(\bar{m}, z)).
\eeq

We define the Lagragian functional
\beq\label{eq:Lagrangian Inner}
\cL(u,v,m,z) = Q(u) + r(u,v,m,z),
\eeq
where the adjoint $v$ is regarded as the Lagrangian multiplier of the state equation.

By requiring the first order variation of \eqref{eq:Lagrangian Inner} at $\bar{m}$ with respect to the adjoint $v$ to vanish, we obtain the state problem: find $u \in \cU$, such that 
\beq\label{eq:state}
\langle \tilde{v}, \partial_v \bar{r}\rangle = 0, \quad \forall \tilde{v} \in \cV.
\eeq
Similarly, by setting the first order variation of \eqref{eq:Lagrangian Inner} at $\bar{m}$ with respect to the state $u$ to zero, 
we obtain the adjoint problem: find $v \in \cV$, such that
\beq\label{eq:adjoint}
\langle \tilde{u} , \partial_u \bar{r}\rangle = - \langle \tilde{u}, \partial_u \bar{Q} \rangle, \quad \forall \tilde{u} \in \cU.
\eeq 

Then, the gradient of $Q$ at $\bar{m}$ acting in the direction $\tilde{m} = m - \bar{m}$ is given by  
\beq\label{eq:gradient}
\langle \tilde{m}, \func{Q}_m(\bar{m}) \rangle  =  \langle \tilde{m}, \partial_m \bar{r}\rangle,
\eeq
where $u$ solves the state problem \eqref{eq:state} and $v$ solves the adjoint problem \eqref{eq:adjoint}.

Similarly, for the computation of the Hessian of $Q$ acting on $\hat{m}$, we consider the Lagrangian functional
\beq\label{eq:meta_lag}
\cL^H(u,v,m,z; \hat{u}, \hat{v}, \hat{m}) = \langle \hat{m}, \partial_m \bar{r}\rangle + \langle \hat{v}, \partial_v \bar{r}\rangle + \langle \hat{u} , \partial_u \bar{r} + \partial_u \bar{Q} \rangle,
\eeq
where $\hat{u}$ and $\hat{v}$ denote the incremental state and incremental adjoint, respectively.
By taking variation of \eqref{eq:meta_lag} with respect to the adjoint $v$ and using \eqref{eq:adjoint_op}, we obtain the incremental state problem: find $\hat{u} \in \cU$ such that
\beq\label{eq:incfwd}
\langle \tilde{v},  \partial_{vu} \bar{r}\, \hat{u} \rangle = - \langle \tilde{v}, \partial_{vm} \bar{r}\, \hat{m}\rangle, \quad \forall \tilde{v} \in \cV,
\eeq 
where the derivatives $\partial_{vu} \bar{r}: \cU \to \cV'$ and $\partial_{vm} \bar{r}: \cM \to \cV'$ are linear operators.
The incremental adjoint problem, obtained by taking variation of \eqref{eq:meta_lag} with respect to the state $u$ and using \eqref{eq:adjoint_op}, reads: find $\hat{v} \in \cV$ such that
\beq\label{eq:incadj}
\langle \tilde{u}, \partial_{uv} \bar{r} \,\hat{v} \rangle  = - \langle \tilde{u}, \partial_{uu} \bar{r}\,\hat{u} + \partial_{uu} \bar{Q}\,\hat{u}  + \partial_{um} \bar{r}\,\hat{m}\rangle, \quad \forall \tilde{u} \in \cU,
\eeq
where $\partial_{uv} \bar{r}: \cV \to \cU'$ is the adjoint of $\partial_{vu} \bar{r}: \cU \to \cV'$ in the sense of \eqref{eq:adjoint_op}.
We remark that an extra term $-\langle \tilde{u}, \partial_{um} \bar{Q} \,\hat{m}\rangle$ is added on the right hand side if $Q$ explicitly depends on $m$.
The Hessian of $Q$ at $\bar{m}$ acting on $\hat{m}$ can then be computed by taking variation of \eqref{eq:meta_lag} with respect to $m$ and using \eqref{eq:adjoint_op} as 
\beq\label{eq:HessianQ}
\langle \tilde{m}, \oper{Q}_{mm}(\bar{m})\,\hat{m}\rangle  = \langle \tilde{m}, \partial_{mv}\bar{r} \, \hat{v} + \partial_{mu}\bar{r}\,\hat{u}  + \partial_{mm} \bar{r}\,\hat{m}\rangle,
\eeq 
where the incremental state $\hat{u}$ and adjoint $\hat{v}$ solve \eqref{eq:incfwd} and \eqref{eq:incadj}, respectively.
Note that $\langle \tilde{m}, \partial_{mu}\bar{Q}\,\hat{u} + \partial_{mm} \bar{Q}\,\hat{m}\rangle$ is also added if $Q$ explicitly depends on $m$.

We summarize the computation of the gradient and Hessian of the control objective $Q$ with respect to the uncertain parameter $m$ in Algorithm \ref{alg:gradientHessian}. 

\begin{algorithm}
\caption{Computation of the gradient and Hessian of $Q$ with respect to $m$.}
\label{alg:gradientHessian}
\begin{algorithmic}
\STATE{\textbf{Input: } control objective $Q$, PDE $r(u,v,m,z)$, mean $\bar{m}$, direction $\hat{m}$.}
\STATE{\textbf{Output: } the gradient $Q_m(\bar{m})$ and Hessian $Q_{mm}(\bar{m})$ acting on $\hat{m}$.}
\vspace*{0.2cm}
\STATE{1. Solve the state problem \eqref{eq:state} for $u$.}
\STATE{2. Solve the adjoint problem \eqref{eq:adjoint} for $v$.}
\STATE{3. Compute the gradient $Q_m(\bar{m})$ by \eqref{eq:gradient}.}
\vspace*{0.2cm}
\STATE{4. Solve the incremental state problem \eqref{eq:incfwd} for $\hat{u}$.}
\STATE{5. Solve the incremental adjoint problem \eqref{eq:incadj} for $\hat{v}$.}
\STATE{6. Compute the Hessian $Q_{mm}(\bar{m})$ acting on $\hat{m}$ by \eqref{eq:HessianQ}.}
\end{algorithmic}
\end{algorithm}

\subsection{Randomized algorithms to compute the trace}
\label{subsec:randtrace}

The first and second moments of the quadratic approximation $Q_{\text{quad}}$ in equation \eqref{eq:quadraticE} involve the trace of the covariance preconditioned Hessian $\cH$. $\cH$ is formally a large scale dense operator and it is only implicitly defined through its action on a given direction as in \eqref{eq:HessianQ}. 

The standard approach to estimate the trace of an implicitly defined operator is based on a Monte Carlo method, where one computes the action of the operator on $N$ random directions sampled from a suitable distribution. For example, using the Gaussian trace estimator, one can approximate $\trace(\cH)$ and $\trace(\cH^2)$ as
\beq\label{eq:GaussTrace1}
\trace(\cH) \approx \widehat{T}_1(\cH) := \frac{1}{N}\sum_{j=1}^{N} \langle \hat{m}_j,  \oper{Q}_{mm}(\bar{m})\,\hat{m}_j \rangle, 
\eeq
and 
\beq\label{eq:GaussTrace}
\trace(\cH^2) \approx \widehat{T}_1(\cH^2) := \frac{1}{N}\sum_{j=1}^{N} \langle  \oper{Q}_{mm}(\bar{m})\,\hat{m}_j, \cC \oper{Q}_{mm}(\bar{m})\,\hat{m}_j\rangle,
\eeq
where $\hat{m}_j$, $j= 1, \dots, N$, are \emph{i.i.d.} random samples from the Gaussian distribution $\cN(0, \cC)$.
Lower bounds on the number of samples $N$ to obtain a guaranteed probabilistic error bound for the trace are discussed in \cite{AvronToledo2011}. Such bounds hold for any symmetric positive defined operator, regardless of size and spectral properties. However the number of samples $N$ necessary to obtain an accurate approximation may be prohibitively large and impractical: e.g., to guarantee a relative error of $10^{-3}$ in the estimation of the trace more than $10^6$ samples are need using the bounds in \cite{AvronToledo2011}. Improved bounds are demonstrated in \cite{Roosta-KhorasaniAscher2015} and show that the Gaussian estimator becomes very inefficient if the stable rank of the operator is small and that it allows for small $N$ only if the eigenvalues are all of  approximately the same size. This means that the randomized Gaussian estimator is not a viable solution to estimate the trace of $\cH$, since it has been observed numerically or proven analytically that for many problems the Hessian operator is either nearly low-rank or its eigenvalues exhibit fast decay \cite{BashirWillcoxGhattasEtAl08, FlathWilcoxAkcelikEtAl11,
  Bui-ThanhGhattas12a, Bui-ThanhGhattas13a, Bui-ThanhGhattas12,
  Bui-ThanhBursteddeGhattasEtAl12_gbfinalist,
  Bui-ThanhGhattasMartinEtAl13, ChenVillaGhattas2017,
  AlexanderianPetraStadlerEtAl16, AlexanderianPetraStadlerEtAl17,
  AlexanderianPetraStadlerEtAl14, CrestelAlexanderianStadlerEtAl17,
  PetraMartinStadlerEtAl14, IsaacPetraStadlerEtAl15,
  MartinWilcoxBursteddeEtAl12, Bui-ThanhGhattas15, ChenGhattas18b}.


In this work, we propose to 
estimate the trace of $\cH$ and $\cH^2$ by means of a randomized approximated eigendecomposition. 
Specifically, we approximate $\trace(\cH)$ and $\trace(\cH^2)$ by 
\beq\label{eq:RandomizedEigen} 
\trace(\cH) \approx \widehat{T}_2(\cH) := \sum_{j=1}^{N} \lambda_j(\cH), \text{ and } \trace(\cH^2) \approx \widehat{T}_2(\cH^2) := \sum_{j=1}^{N} \lambda_j^2(\cH),
\eeq
where $\lambda_j(\cH)$, $j = 1, \dots, N$, are the dominant eigenvalues of $\cH$ obtained by solving the generalized eigenvalue problem
\beq\label{eq:gEigen}
\langle \phi, \oper{Q}_{mm}(\bar{m}) \, \psi_j \rangle = \lambda_j \langle \phi, \cC^{-1} \psi_j \rangle, \quad \forall \phi \in \cM, \quad j = 1, \dots, N.
\eeq
Above, the eigenvector $\psi_j$, $j = 1, \dots, N$, are $\cC^{-1}$-orthonormal, that is
\beq
\langle \psi_j, \cC^{-1} \psi_i\rangle = \delta_{ij}, \quad i,j = 1, \dots, N,
\eeq 
where $\delta_{ij}$ is the Kronecker delta. Note that $N$ is independent of the nominal parameter dimension, and, as shown numerically in Sec.\ \ref{sec:numerics}, it is relatively small for $\cH$ nearly low-rank or when its spectrum exhibits fast decay.

After suitable discretization, e.g., by a finite element discretization, the generalized eigenvalue problem \eqref{eq:gEigen} results in an algebraic eigenproblem of the form $A \psi = \lambda B \psi$ with $A, B \in \bR^{n\times n }$ and $\psi  \in \bR^n$.
Algorithm \ref{alg:randomizedEigenSolver} summarizes the so-called \emph{double pass randomized algorithm} to solve the algebraic generalized eigenvalue problem (see \cite{HalkoMartinssonTropp2011, SaibabaLeeKitanidis2016} for details of the algorithms and \cite{VillaPetraGhattas2017} for its implementation). 

\begin{algorithm}
\caption{Randomized algorithm for generalized eigenvalue problem $(A, B)$}
\label{alg:randomizedEigenSolver}
\begin{algorithmic}
\STATE{\textbf{Input: } matrix $A, B$, the number of eigenpairs $k$, an oversampling factor $p$.}
\STATE{\textbf{Output: } eigenpairs $(\Lambda_k, \Psi_k)$ with $\Lambda_k = \text{diag}(\lambda_1, \dots, \lambda_k)$ and $\Psi_k = (\psi_1, \dots, \psi_k)$.}
\vspace*{0.2cm}
\STATE{1. Draw a Gaussian random matrix $\Omega \in \bR^{n\times (k+p)}$.}
\STATE{2. Compute $Y = B^{-1} (A \Omega)$.}
\STATE{3. Compute $QR$-factorization $Y = QR$ such that $Q^\top B Q = I_{k+p}$.}
\STATE{4. Form $T = Q^\top A Q$ and compute eigendecomposition $T = S \Lambda S^\top$.}
\STATE{5. Extract $\Lambda_k = \Lambda(1:k, 1:k)$ and $\Psi_k = QS_k$ with $S_k = S(:,1:k)$.}
\end{algorithmic}
\end{algorithm}

{
There are multiple advantages in using Algorithm \ref{alg:randomizedEigenSolver} to approximate the trace, see e.g., \cite{HalkoMartinssonTropp2011, SaibabaAlexanderianIpsen2016}. In terms of accuracy, the approximation error is bounded by the sum of the remaining eigenvalues, so that the error is small if the eigenvalues decay fast or if the Hessian $\cH$ has low rank, see \cite{HalkoMartinssonTropp2011, SaibabaAlexanderianIpsen2016} for more details. 
In terms of computational efficiency, the $2(k+p)$ Hessian matrix-vector products, which entail the solution of a pair of linearized state/adjoint equations (as shown in Sec.\ \ref{subsec:gradientHessian}) and therefore dominate the computational cost of the algorithm, can be computed independently and, therefore, asynchronously across the random directions.}
\subsection{Monte Carlo correction}
\label{sec:MCVR}
Using the moments of the truncated Taylor expansion $Q_{\text{lin}}$ and $Q_{\text{quad}}$ as approximation of $\bE[Q]$ and $\Var[Q]$ in the cost functional \eqref{eq:objective} 
may be not sufficiently accurate.
In this circumstance, we can use the Taylor approximation as a control variate to reduce the variance of the Monte Carlo estimator. The basic idea is to correct the moments of $Q_{\text{lin}}$ and $Q_{\text{quad}}$ by applying a Monte Carlo method to estimate the mean (and variance) of the difference between $Q$ and its Taylor approximation.

Specifically, recalling \eqref{eq:linearE} and \eqref{eq:linearQ}, the Monte Carlo correction for the mean of the linear approximation reads
\beq\label{eq:linearMC}
\begin{split}
\bE[Q]  &= \bE[Q_{\text{lin}}] + \bE[Q - Q_{\text{lin}}] 
 \\
 &\approx \widehat{Q}_{\rm lin} := Q(\bar{m}) + \frac{1}{M_1} \sum_{i = 1}^{M_1} \Big(Q(m_i) - Q(\bar{m}) - \langle m_i - \bar{m}, \func{Q}_m(\bar{m})\rangle\Big).
\end{split}
\eeq
Similarly, recalling that the variance is invariant with respect to deterministic translations, we have
\beq
\begin{aligned}
\Var[Q] 
&= \Var[Q-Q(\bar{m})] 
= \bE[(Q-Q(\bar{m}))^2] - (\bE[Q-Q(\bar{m})])^2 
\\
& 
= \bE[(Q_{\text{lin}}-Q(\bar{m}))^2] + \bE[(Q-Q(\bar{m}))^2 - (Q_{\text{lin}}-Q(\bar{m}))^2] 
\\
& 
- \big(\bE[Q_{\text{lin}}-Q(\bar{m})] + \bE[(Q-Q(\bar{m})) - (Q_{\text{lin}}-Q(\bar{m}))]\big)^2,
\end{aligned}
\eeq
and therefore the Monte Carlo correction of the variance of the linear approximation is given by
\beq\label{eq:linMCVar}
\begin{aligned}
\widehat{V}_Q^\text{lin} 
& 
:= \langle \cC \func{Q}_m(\bar{m}), \func{Q}_m(\bar{m}) \rangle 
+ \frac{1}{M_2} \sum_{i=1}^{M_2} \Big((Q(m_i)-Q(\bar{m}))^2 - (\langle m_i - \bar{m}, \func{Q}_m(\bar{m})\rangle)^2\Big) 
\\
& 
- \left(\frac{1}{M_2} \sum_{i=1}^{M_2} Q(m_i) - Q(\bar{m}) - \langle m_i - \bar{m}, \func{Q}_m(\bar{m})\rangle\right)^2.
\end{aligned}
\eeq
Note that, similar to the sample average approximation in Sec.\ \ref{sec:MCinteg}, for simplicity we can use the same $M = M_1 = M_2 $ random samples for the approximation of both the mean and  variance. A more careful choice of $M_1$ and $M_2$ relies on the variance of $Q^{(1)} = Q - Q_{\text{lin}}$ and $Q^{(2)} = (Q-\bar{Q})^2 - (Q-Q_{\text{lin}})^2$, 
in order to balance the errors in the approximation of the mean and variance.
For the quadratic approximation, a combination of \eqref{eq:quadraticE} and \eqref{eq:quadraticQ} leads to 
\beq
\begin{aligned}\label{eq:quadraticMC}
\bE[Q]  &  = \bE[Q_{\text{quad}}] + \bE[Q - Q_{\text{quad}}]   \approx \widehat{Q}_{\rm quad} := Q(\bar{m}) + \frac{1}{2}\text{tr}(\cH) \\
& + \frac{1}{M_1} \sum_{i = 1}^{M_1}  \Big(Q(m_i) - Q(\bar{m}) - \langle m_i - \bar{m}, \func{Q}_m(\bar{m}) \rangle -  \frac{1}{2}\langle  m_i - \bar{m}, \oper{Q}_{mm}(\bar{m})\,(m_i - \bar{m}) \rangle \Big).
\end{aligned}
\eeq
For the variance, akin to the linear approximation, we obtain 
\beq
\begin{aligned}
\Var[Q] & = \bE[(Q_{\text{quad}}-Q(\bar{m}))^2] + \bE[(Q-Q(\bar{m}))^2 - (Q_{\text{quad}}-Q(\bar{m}))^2] \\
& - \big(\bE[Q_{\text{quad}}-Q(\bar{m})] + \bE[(Q-Q(\bar{m})) - (Q_{\text{quad}}-Q(\bar{m}))]\big)^2,\\
\end{aligned}
\eeq
which can be approximated by
\beq\label{eq:quadMCVar}
\begin{aligned}
\widehat{V}_Q^\text{quad} & := \langle \cC \func{Q}_m(\bar{m}), \func{Q}_m(\bar{m}) \rangle + \frac{1}{4} (\text{tr}(\cH))^2 + \frac{1}{2} \text{tr}(\cH^2)\\
& + \frac{1}{M_2} \sum_{i=1}^{M_2} \Big((Q(m_i) -Q(\bar{m}))^2  -   \Big( \langle m_i - \bar{m}, \func{Q}_m(\bar{m}) \rangle + \frac{1}{2} \langle m_i - \bar{m} ,\oper{Q}_{mm}(\bar{m})\,(m_i - \bar{m}) \rangle \Big)^2\Big)\\
& - \Big(
\frac{1}{2}\text{tr}(\cH) + \frac{1}{M_2}\sum_{i=1}^{M_2} \Big(Q(m_i) - \langle m_i - \bar{m}, \func{Q}_m(\bar{m}) \rangle    - \frac{1}{2}\langle m_i - \bar{m} ,\oper{Q}_{mm}(\bar{m})\,(m_i - \bar{m}) \rangle \Big)
\Big)^2.
\end{aligned}
\eeq

Note that to achieve certain approximation accuracy, the number of Monte Carlo samples depends on the variance of the integrand.
If the approximations $Q_{\text{lin}}$ and $Q_{\text{quad}}$ are highly correlated with $Q$, the variances of $Q - Q_{\text{lin}}$ and $Q - Q_{\text{quad}}$ would be small so that only a small number of Monte Carlo samples are required. 
This use of the Taylor approximation (of the control objective with respect to the uncertain parameter) can be viewed as a kind of ``multifidelity Monte Carlo method" \cite{PeherstorferWillcoxGunzburger18}, in which the Taylor approximation plays the role of a reduced model.
\section{Gradient-based optimization}
\label{sec:gradient}
\label{sec:gradientopt}
In this section, we present a gradient-based optimization method to solve the problem \eqref{eq:optimization} using different approximations introduced in Sec.\ \ref{sec:taylor}, namely the sample average approximation, the first and second order Taylor approximation, and variance reduction using the Taylor approximation as the control variate. These approximated formulations feature different levels of accuracy and computational complexity, and could possibly be combined within a multifidelity framework to obtain an efficient computational procedure while preserving high accuracy. 
In this section, we use the Lagrangian formalism to derive expressions for the gradient of the approximated formulations of the cost functional with respect to the control $z$, which are needed to perform the gradient-based optimization.
In the following, to distinguish the gradient of the cost function with respect to the control from the gradient of the control objective with respect to the uncertain parameter, we denote the former as $z$-gradient.

\subsection{$z$-gradient for the sample average approximation}
By the sample average approximation presented in Sec.\ \ref{sec:MCinteg}, the cost functional  becomes 
\beq\label{eq:objectiveMC}
\cJ^{\text{MC}}(z) = \widehat{Q} + \beta \widehat{V}_Q + \cP(z),
\eeq
subject to the state problem \eqref{eq:PDE} at $m_i$, $i = 1, \dots, M$, where $M = M_1 = M_2$ if the same samples are used for the mean and variance or $M = M_1 + M_2$ otherwise. 
To compute the gradient of \eqref{eq:objectiveMC} with respect to the control $z$ ($z$-gradient for short), we consider the Lagrangian functional
\beq
\cL^{\text{MC}}(\{u_i\},\{v_i\},z) = \cJ^{\text{MC}}(z) + \sum_{i=1}^M r(u_i,v_i,m_i,z),
\eeq
where $\{v_i\}$ are the Lagrange multipliers (adjoints) of the state problems: given $m_i \in \cM$ and $z\in \cZ$, find $u_i \in \cU$ such that
\beq \label{eq:forwardMC}
r(u_i, v_i, m_i,z) = 0, \quad \forall v_i \in \cV; \quad i = 1, \ldots M.
\eeq
The $i$-th adjoint problem ($i=1, \ldots M$) is obtained by setting to zero the first variation of the Lagrangian functional with respect to the state $u_i$. The $i$-th adjoint $v_i$ then solves: given $m_i \in \cM$, $z\in \cZ$, and $u_i \in \cU$ as the solution of \eqref{eq:forwardMC}, find $v_i \in \cV$ such that
\beq \label{eq:adjointMC}
\begin{aligned}
\langle \tilde{u}_i, \partial_{u} r(u_i,v_i,m_i,z) \rangle & = - \frac{1}{M_1} \langle \tilde{u}_i, \partial_{u} Q(m_i) \rangle - \frac{2\beta}{M_2} Q(m_i)\langle \tilde{u}_i, \partial_u Q(m_i)\rangle \\
& + \frac{2\beta}{M_2} \left(\frac{1}{M_2} \sum_{j=1}^{M_2} Q(m_j) \right) \langle \tilde{u}_i, \partial_u Q(m_i) \rangle, \quad \forall \tilde{u}_i \in \cU.
\end{aligned}
\eeq
Then we can compute the $z$-gradient as 
\beq
\begin{aligned}
\langle \tilde{z}, D_z \cJ^{\text{MC}}(z) \rangle & = \langle \tilde{z}, \partial_z \cL^{\text{MC}}(\{u_i\}, \{v_i\}, z) \rangle  = \langle\tilde{z}, D_z \cP(z)\rangle + \sum_{i = 1}^M \langle \tilde{z}, \partial_z r(u_i, v_i, m_i, z)\rangle,
\end{aligned}
\eeq
where $\{u_i\}$ and $\{v_i\}$ denote the set of solutions of the state problem \eqref{eq:forwardMC} and adjoint problem \eqref{eq:adjointMC} at the set of samples $\{m_j\}$, respectively. Therefore, each gradient evaluation requires $M$ (nonlinear) state PDE solves and $M$ linearized PDE solves.
{Note that, for the general case in which $Q$ explicitly depends on $z$, extra terms involving $\langle \tilde{z}, \partial_z Q(m_i) \rangle$ need to be included in the gradient computation.}

\subsection{$z$-gradient for the Taylor approximation}\label{subsec:GradientTaylor}
We first consider the linear approximation \eqref{eq:linearE}, by which the cost functional \eqref{eq:objective} becomes   
\beq\label{eq:Jlin}
\cJ_{\text{lin}}(z) = \bar{Q} + \beta \langle \cC \func{Q}_m(\bar{m}), \func{Q}_m(\bar{m}) \rangle + \cP(z),
\eeq
where $\bar{Q}$ and its gradient $\func{Q}_m(\bar{m})$ are defined in \eqref{eq:r0Q0} and \eqref{eq:gradient}, respectively.

To compute the $z$-gradient of \eqref{eq:Jlin}, we introduce the Lagrangian functional 
\beq\label{eq:LagrangianLinear}
\begin{aligned}
\cL_{\text{lin}}(u,v,u^*, v^*, z) &= \cJ_{\text{lin}}(z) + \langle v^*, \partial_v \bar{r} \rangle + \langle  u^*, \partial_u \bar{r} + \partial_u \bar{Q} \rangle,
\end{aligned}
\eeq
where the adjoints $v^* \in \cV$ and $u^*\in \cU$ represent the Lagrange multipliers of the state problem \eqref{eq:state} and adjoint problem \eqref{eq:adjoint}.

By setting to zero the first order variation of the Lagrangian \eqref{eq:LagrangianLinear} with respect to $v$ and using \eqref{eq:adjoint_op}, we have: find $u^* \in \cU$ such that 
\beq\label{eq:adjointLinear}
\langle \tilde{v} , \partial_{vu}\bar{r}  u^* \,\rangle =  -2\beta  \langle \tilde{v}, \partial_{vm} \bar{r} \, (\cC \partial_m \bar{r})  \rangle, \quad \forall \tilde{v} \in \cV.
\eeq
while setting to zero the first order variation with respect to $u$ and using \eqref{eq:adjoint_op}, we have: find $v^* \in \cV$ such that 
\beq\label{eq:adjointLinear2}
\begin{aligned}
\langle \tilde{u}, \partial_{uv}\bar{r} \, v^* \rangle = & -\langle \tilde{u}, \partial_u \bar{Q} \rangle - 2\beta \langle \tilde{u} , \partial_{um} \bar{r} \, (\cC \partial_m \bar{r}) \rangle  - \langle \tilde{u}, \partial_{uu}\bar{r} \, u^*  + \partial_{uu}\bar{Q} \, u^* \rangle, \quad \forall \tilde{u} \in \cU.
 \end{aligned}
\eeq
Finally, the $z$-gradient of $\cJ_{\text{lin}}$ in a direction $\tilde{z} \in \cZ$ is computed as
\beq\label{eq:linearGradient}
\begin{aligned}
\langle \tilde{z}, D_z \cJ_{\text{lin}}(z) \rangle & = \langle \tilde{z}, \partial_z \cL_{\text{lin}}(u,v,u^*, v^*, z)\rangle  \\
& = 2 \beta  \langle \tilde{z}, \partial_{zm} \bar{r} \,(\cC \partial_m \bar{r}) \rangle + \langle \tilde{z}, \partial_z \cP(z) \rangle  +  \langle \tilde{z}, \partial_{zv} \bar{r} \, v^* + \partial_{zu} \bar{r} \, u^* \rangle.
\end{aligned}
\eeq
In summary, evaluation of the cost functional with the linear approximation requires solution of the state problem \eqref{eq:state} and adjoint problem \eqref{eq:adjoint}. Computation of the $z$-gradient requires, in addition, solution of the two linear problems in \eqref{eq:adjointLinear} and \eqref{eq:adjointLinear2}.

For the quadratic approximation \eqref{eq:quadraticQ}, thanks to the relation \eqref{eq:quadraticE} and the approximation \eqref{eq:GaussTrace}, we obtain the approximation of the cost functional as
\beq\label{eq:Jquad}
\cJ_{\text{quad}}(z) = \bar{Q} + \frac{1}{2} \sum_{j=1}^N \lambda_j + \beta \left(\langle \cC \func{Q}_m(\bar{m}), \func{Q}_m(\bar{m}) \rangle + \frac{1}{2} \sum_{j=1}^N \lambda_j^2 \right) +\cP(z),
\eeq
which is subject to the state and adjoint problems \eqref{eq:state} and \eqref{eq:adjoint} for the evaluation of $Q$ and its gradient at $\bar{m}$, the generalized eigenvalue problems \eqref{eq:gEigen} for the computation of the eigenvalues, as well as the incremental state and adjoint problems \eqref{eq:incfwd} and \eqref{eq:incadj} for the Hessian action in \eqref{eq:gEigen}. Correspondingly, we form the Lagrangian 
\beq\label{eq:LagrangianQuadratic}
\begin{aligned}
& \cL_{\text{quad}} \left(u,v, \{\lambda_j\}, \{\psi_j\}, \{\hat{u}_j\}, \{\hat{v}_j\}, u^*,v^*, \{\lambda_j^*\}, \{\psi_j^*\}, \{\hat{u}^*_j\}, \{\hat{v}^*_j\}, z\right)  \\
&
= \cJ_{\text{quad}}(z)  + \langle v^*,  \partial_v \bar{r} \rangle  + \langle u^*, \partial_u \bar{r} + \partial_u \bar{Q}\rangle   \\
& + \sum_{j=1}^N \langle \psi_j^*, (\oper{Q}_{mm}(\bar{m}) - \lambda_j \cC^{-1}) \psi_j \rangle + \sum_{j = 1}^N \lambda_j^* (\langle \psi_j, \cC^{-1} \psi_j \rangle - 1) \\ 
& +  \sum_{j=1}^{N} \langle \hat{v}^*_j, \partial_{vu} \bar{r} \, \hat{u}_j + \partial_{vm} \bar{r}\, \psi_j\rangle  + \sum_{j=1}^{N} \langle \hat{u}^*_j, \partial_{uv} \bar{r} \,\hat{v}_j + \partial_{uu} \bar{r} \,\hat{u}_j + \partial_{uu} \bar{Q} \,\hat{u}_j +\partial_{um}  \bar{r}\,\psi_j \rangle.
\end{aligned}
\eeq
Here we assume that the dominating eigenvalues are not repeated (thanks to their rapid decay as shown in Sec.\ \ref{sec:numerics}), so that the constraints $\langle \psi_j, \cC^{-1} \psi_j \rangle = 1$, $j = 1, \dots, N$, are sufficient to guarantee the orthonormality $\langle \psi_j, \cC^{-1} \psi_i \rangle = \delta_{ij}$ for any $i,j = 1, \dots, N$. In fact, for any $i \neq j$, by definition we have 
\beq
\langle \psi_j, Q_{mm} \psi_i \rangle = \langle \psi_j, \lambda_i C^{-1} \psi_i \rangle \text{ and } \langle \psi_i, Q_{mm} \psi_j \rangle = \langle \psi_i, \lambda_j \cC^{-1} \psi_j \rangle, 
\eeq 
which, by the symmetry of $Q_{mm}$ and $\cC^{-1}$, leads to 
\beq
(\lambda_i - \lambda_j) \langle \psi_i, \cC^{-1} \psi_j \rangle = 0, \text{ so that } \langle \psi_i, \cC^{-1} \psi_j \rangle = 0 \text{ if }\lambda_i \neq \lambda_j.
\eeq

By setting the variation of $\cL_{\text{quad}}$ with respect to $\lambda_j$ as zero, we obtain 
\beq\label{eq:psijhat}
\psi_j^* = \frac{1+2\beta \lambda_j}{2} \psi_j, \quad j = 1, \dots, N.
\eeq
Subsequently, for each $j = 1, \dots, N$, taking the variation of $\cL_{\text{quad}}$ with respect to $\hat{v}_j$ as zero, and using the Hessian action equation \eqref{eq:HessianQ} and \eqref{eq:adjoint_op}, we have: find $\hat{u}^*_j \in \cU$ such that 
\beq\label{eq:uhatj}
\langle \tilde{v}, \partial_{vu} \bar{r} \, \hat{u}^*_j \rangle =  - \langle \tilde{v}, \partial_{vm} \bar{r} \, \psi^*_j \rangle, \quad \forall \tilde{v} \in \cV,
\eeq
which is the same as the incremental forward problem \eqref{eq:incfwd}. Therefore, from the result \eqref{eq:psijhat}, we have 
\beq\label{eq:hatujstar}
\hat{u}^*_j  = \frac{1+2\beta \lambda_j}{2} \hat{u}_j, \quad j = 1, \dots, N\;.
\eeq
Then, taking variation with respect to $\hat{u}_j$ as zero and using \eqref{eq:adjoint_op}, we have: find $\hat{v}_j^* \in \cV$ such that 
\beq\label{eq:vhatj}
\langle \tilde{u}, \partial_{uv} \bar{r} \, \hat{v}^*_j \rangle = - \langle \tilde{u}, \partial_{uu} \bar{r} \, \hat{u}^*_j + \partial_{uu} \bar{Q} \, \hat{u}^*_j + \partial_{um} \bar{r} \, \psi^*_j\rangle, \quad \forall \tilde{u} \in \cU,
\eeq 
which is the same as the incremental adjoint problem \eqref{eq:incadj}. Therefore, by \eqref{eq:hatujstar} and \eqref{eq:psijhat}, we obtain 
\beq\label{eq:hatvjstar}
\hat{v}_j^* =  \frac{1+2\beta \lambda_j}{2} \hat{v}_j, \quad j = 1, \dots, N.
\eeq
By setting to zero the variation with respect to $v$ and using \eqref{eq:adjoint_op}, we have: find $u^* \in \cU$ such that 
\beq\label{eq:uhat}
\begin{aligned}
\langle \tilde{v}, \partial_{vu} \bar{r} \, u^* \rangle & = - 2\beta \langle \tilde{v}, \partial_{vm} \bar{r} \,(\cC \partial_m \bar{r}) \rangle 
 -  \sum_{j = 1}^N \langle \tilde{v}, \partial_{vmu} \bar{r} \,  \hat{u}_j \, \psi_j^*+ \partial_{vmm} \bar{r}  \, \psi_j \, \psi_j^* \rangle \\
& - \sum_{j = 1}^N \langle \tilde{v}, \partial_{vuu} \bar{r} \, \hat{u}_j \, \hat{u}_j^* + \partial_{vum} \bar{r} \, \psi_j \, \hat{u}_j^* \rangle, \quad \forall \tilde{v} \in \cV.
\end{aligned}
\eeq
Lastly, by setting to zero the variation with respect to $u$ and using \eqref{eq:adjoint_op}, we have: find $v^* \in \cV$ such that 
\beq\label{eq:vhat}
\begin{aligned}
 \langle \tilde{u}, \partial_{uv} \bar{r} \, v^* \rangle  &=  - \langle \tilde{u}, \partial_u \bar{Q} \rangle - 2\beta \langle \tilde{u}, \partial_{um} \bar{r} \,(\cC \partial_m \bar{r}) \rangle  - \langle \tilde{u}, \partial_{uu} \bar{r} \, u^* + \partial_{uu} \bar{Q} \, u^* \rangle \\
& - \sum_{j = 1}^N \langle \tilde{u}, \partial_{umv} \bar{r} \, \hat{v}_j \, \psi_j^* + \partial_{umu} \bar{r} \, \hat{u}_j \, \psi_j^* + \partial_{uum} \bar{r} \, \psi_j \, \psi_j^* \rangle  \\
& - \sum_{j = 1}^N \langle \tilde{u}, \partial_{uvu} \bar{r} \, \hat{u}_j \, \hat{v}_j^* + \partial_{uvm} \bar{r} \, \psi_j \, \hat{v}_j^*\rangle \\
& - \sum_{j = 1}^N \langle \tilde{u}, \partial_{uuv} \bar{r} \,\hat{v}_j \, \hat{u}_j^* + \partial_{uuu} \bar{r} \,\hat{u}_j \, \hat{u}_j^*+ \partial_{uuu} \bar{Q} \,\hat{u}_j \, \hat{u}_j^*+\partial_{uum}  \bar{r}\,\psi_j\, \hat{u}_j^* \rangle, \quad \forall \tilde{u} \in \cU,
\end{aligned}
\eeq
where the third order derivatives take the form of the bilinear map $\partial_{xyz} f: \cZ \times \cY \to \cX'$.
Then, the $z$-gradient of the cost functional in the direction $\tilde{z}\in \cZ$ can be computed as
\beq\label{eq:DzJquad}
\langle \tilde{z}, D_z \cJ_{\text{quad}}(z)\rangle =  \langle \tilde{z}, \partial_z \cL_{\text{quad}} \left(\text{states}, \text{adjoints}, z\right) \rangle,
\eeq 
where states and adjoints represent the state and adjoint variables for short. Note that we do not need to solve for $\lambda_j^*$ as the terms involving $\lambda_j^*$ are independent of $z$.
In summary, to get these states and adjoints, we need to solve $1$ (nonlinear) state problem and $2N + 3$ linear problems ($1$ for $v$, $N$ for $\{\hat{u}_j\}$, $N$ for $\{\hat{v}_j\}$, 
$1$ for $u^*$, and $1$ for $v^*$), where the linear operators are the same as those in \eqref{eq:incfwd} and \eqref{eq:incadj}. Moreover, to compute the eigenvalues by Algorithm \ref{alg:randomizedEigenSolver}, which requires the action of the Hessian on $N+p$ random directions twice, we need to solve $4(N+p)$ linear problems.

We summarize the evaluation of the cost functional and its $z$-gradient with quadratic approximation as follows. 

\begin{algorithm}
\caption{Computation of the cost functional $\cJ_{\text{quad}}$ and its $z$-gradient.}
\label{alg:costGradientQuad}
\begin{algorithmic}
\STATE{\textbf{Input: } cost functional $\cJ(z)$ and PDE $r(u,v,m,z)$.}
\STATE{\textbf{Output: } the approximate cost functional $\cJ_{\text{quad}}$ and its $z$-gradient $D_z \cJ_{\text{quad}}$.}
\vspace*{0.2cm}
\STATE{1. Solve the state problem \eqref{eq:state} for $u$ at $\bar{m}$ and evaluate $\bar{Q} = Q(u(\bar{m}))$.}
\STATE{2. Solve the adjoint problem \eqref{eq:adjoint} for $v$ and evaluate $Q_m(\bar{m})$ by \eqref{eq:gradient}.}
\STATE{3. Use Algorithm \ref{alg:gradientHessian} and \ref{alg:randomizedEigenSolver} to compute the generalized eigenpairs $(\lambda_j, \psi_j)_{j=1}^N$.}
\STATE{4. Compute the approximate cost functional $\cJ_{\text{quad}}$ by \eqref{eq:Jquad}.}
\vspace*{0.2cm}
\STATE{5. Solve the linear problem \eqref{eq:uhat} for $\hat{u}^*$.}
\STATE{6. Solve the linear problem \eqref{eq:vhat} for $\hat{v}^*$.}
\STATE{7. Compute the $z$-gradient $D_z \cJ_{\text{quad}}$ by \eqref{eq:DzJquad}.}
\end{algorithmic}
\end{algorithm}

\subsection{$z$-gradient for the Monte Carlo correction}
The approximation in \eqref{eq:Jlin} and \eqref{eq:Jquad} introduced by the Taylor approximation may be not sufficiently accurate and generates bias in the evaluation of the mean $\bE[Q]$ and variance $\text{Var}[Q]$.
Recalling \eqref{eq:linearMC} and \eqref{eq:linMCVar}, we obtain a MC-corrected unbiased linear approximation for the cost functional, which is given by  
\beq\label{eq:JlinMC}
\cJ_{\text{lin}}^{\text{MC}}(z) =  \widehat{Q}_{\text{lin}} + \beta \widehat{V}_Q^{\text{lin}} + \cP(z), 
\eeq
subject to the state problem \eqref{eq:state} at $\bar{m}$, the adjoint problem \eqref{eq:adjoint}, and the state problems \eqref{eq:PDE} at $m_i$, $i = 1, \dots, M$, where $M = M_1 = M_2$ if the same samples are used for the mean and variance or $M = M_1 + M_2$ otherwise.
We form the corresponding Lagrangian functional as 
\beq
\begin{aligned}
 \cL_{\text{lin}}^{\text{MC}}\Big(u, v, \{u_i\}, u^*, v^*, \{v_i\}, z \Big)  = \cJ_{\text{lin}}^{\text{MC}}(z)  + \langle v^*, \partial_v \bar{r} \rangle +  \langle u^*, \partial_u \bar{r} + \partial_u \bar{Q} \rangle   + \sum_{i=1}^{M}  r(u_i, v_i, m_i, z).
\end{aligned} 
\eeq
Setting variation with respect to $v^*$, $v_i$, and $u^*$ to zero we obtain $u$, $u_i$, and $v$ by solving $1+M$ state problems \eqref{eq:PDE} at $m_i$, $i = 0, \dots, M$, and $1$ linear adjoint problem \eqref{eq:adjoint}. 
In a similar way, setting variation with respect to $v$, $u_i$, and $u$ to zero, we obtain $u^*$, $v_i$, and $v^*$, by solving $2+M$ linear problems. 
Thus, in addition, $2+M$ linear problems need to be solved to compute the $z$-gradient
\beq
\langle \tilde{z}, D_z \cJ_{\text{lin}}^{\text{MC}}(z) \rangle = \langle \tilde{z}, \partial_z \cL_{\text{lin}}^{\text{MC}}\Big(u, v, \{u_i\}, u^*, v^*, \{v_i\}, z \Big)\rangle.
\eeq

As for the quadratic approximation \eqref{eq:quadraticQ}, by the Monte Carlo correction \eqref{eq:quadraticMC} and \eqref{eq:quadMCVar}, we have an unbiased and more accurate evaluation of the cost functional as 
\beq\label{eq:JquadMC}
\cJ_{\text{quad}}^{\text{MC}}(z) =  \widehat{Q}_{\text{quad}} + \beta \widehat{V}_Q^{\text{quad}} + \cP(z), 
\eeq
subject to the state problems \eqref{eq:PDE} and \eqref{eq:state}, the adjoint problem \eqref{eq:adjoint}, the incremental state and adjoint problems \eqref{eq:incfwd} and \eqref{eq:incadj}. 
We form the Lagrangian with the Lagrange multiplier $(u^*,v^*,\{v_i\}, \{\lambda^*_j\}, \{\psi^*_j\}, \{\hat{u}^*_j\}, \{\hat{v}^*_j\}, \{\hat{u}^*_i\}, \{\hat{v}^*_i\})$ as 
\beq\label{eq:LagrangianQuadMC}
\begin{aligned}
& \cL_{\text{quad}}^{\text{MC}}\Big(u,v,\{u_i\}, \{\lambda_j\}, \{\psi_j\}, \{\hat{u}_j\}, \{\hat{v}_j\}, \{\hat{u}_i\}, \{\hat{v}_i\},  u^*,v^*,\{v_i\}, \{\lambda^*_j\}, \{\psi^*_j\}, \{\hat{u}^*_j\}, \{\hat{v}^*_j\}, \{\hat{u}^*_i\}, \{\hat{v}^*_i\}, z\Big) \\ 
& = \cJ_{\text{quad}}^{\text{MC}} + \langle v^*, \partial_v \bar{r} \rangle +  \langle u^*, \partial_u \bar{r} + \partial_u \bar{Q} \rangle   + \sum_{i=1}^{M}  r(u_i, v_i, m_i, z) \\
& + \sum_{j=1}^N \langle \psi_j^*, (\oper{Q}_{mm}(\bar{m}) - \lambda_j \cC^{-1}) \psi_j \rangle + \sum_{j = 1}^N \lambda_j^* (\langle \psi_j, \cC^{-1} \psi_j \rangle - 1) \\ 
& +  \sum_{j=1}^{N} \langle \hat{v}^*_j, \partial_{vu} \bar{r} \, \hat{u}_j + \partial_{vm} \bar{r}\, \psi_j\rangle  + \sum_{j=1}^{N} \langle \hat{u}^*_j, \partial_{uv} \bar{r} \,\hat{v}_j + \partial_{uu} \bar{r} \,\hat{u}_j + \partial_{uu} \bar{Q} \,\hat{u}_j +\partial_{um}  \bar{r}\,\psi_j \rangle\\
& +  \sum_{i=1}^{M} \langle \hat{v}^*_i, \partial_{vu} \bar{r} \, \hat{u}_i + \partial_{vm} \bar{r}\, m_i\rangle  + \sum_{i=1}^{M} \langle \hat{u}^*_i, \partial_{uv} \bar{r} \,\hat{v}_i + \partial_{uu} \bar{r} \,\hat{u}_i + \partial_{uu} \bar{Q} \,\hat{u}_i +\partial_{um}  \bar{r}\,m_i \rangle.
\end{aligned}
\eeq
By setting variation with respect to $v^*, u^*, v_i, \hat{v}^*_j, \hat{u}^*_j, \hat{v}^*_i, \hat{u}^*_i$ to zero, we obtain the states $u,v,u_i, \hat{u}_j, \hat{v}_j, \hat{u}_i, \hat{v}_i$ by solving $1+M$ (nonlinear) state problems and $1+2N + 2M$ linear problems, $i = 1, \dots, M$, $j = 1, \dots, N$. For the eigenpairs, $4(N+p)$ linear problems needs to be solved. 
By setting variation with respect to $\lambda_j$, $\hat{v}_j$, and $\hat{u}_j$ to zero, we obtain similar expressions \eqref{eq:psijhat}, \eqref{eq:hatujstar}, and \eqref{eq:hatvjstar} for the adjoint variables $\psi_j^*$, $\hat{u}_j^*$, and $\hat{v}_j^*$, where the constant $\frac{1+2\beta \lambda_j}{2}$ is replaced by $\frac{1+2\beta \lambda_j}{2} - \beta c $ where $c$ is given by
\beq
c = \frac{1}{M} \sum_{i = 1}^M Q(m_i) - Q(\bar{m}) - \langle m_i - \bar{m}, Q_m(\bar{m}) \rangle - \frac{1}{2} \langle m_i - \bar{m}, Q_{mm}(\bar{m})(m_i - \bar{m}) \rangle.
\eeq
 By setting variation with respect to $\hat{v}_i, \hat{u}_i, u_i, v, u$ to vanish, we get the adjoints $\hat{u}^*_i, \hat{v}^*_i, v_i, u^*, v^*$ by solving $2+3M$ linear problems, which are defined similarly as that for the quadratic approximation, but with the additional correction terms in both the cost functional and the constraints, thus in total $1+M$ (nonlinear) state problems and $3+6N+4p + 5M$ linear problems to evaluate the $z$-gradient 
\beq
\langle \tilde{z}, D_z \cJ_{\text{quad}}^{\text{MC}}(z) \rangle = \langle \tilde{z}, \partial_z \cL_{\text{quad}}^{\text{MC}}(\text{states}, \text{adjoints}, z) \rangle,
\eeq
where states and adjoints represent the state and adjoints for short.
\subsection{$z$-gradient-based optimization method} 
\label{subsec:gradientopt}
To solve the PDE-constrained optimal control problem we use a gradient-based optimization method, since it only requires the ability to evaluate the cost functional \eqref{eq:objective} and its $z$-gradient. In particular, we use a BFGS method \cite{NocedalWright2006}.

The five approximations described in this section all differ among each other in terms of accuracy and computational cost. Table \ref{table:Cost} summarizes the cost, measured in number of (nonlinear) state and linear PDE solves, to evaluate cost functional and its $z$-gradient for each approximation. The linear approximation in \eqref{eq:Jlin} is the cheapest to compute but also the most crude, while the MC-corrected approximations in \eqref{eq:JlinMC} and \eqref{eq:JquadMC} are the most accurate but also the most expensive. The quadratic approximation \eqref{eq:Jquad} is, in many application, more accurate then the linear approximation and, for highly nonlinear problems, is much cheaper to compute than the MC-corrected approximation. In fact, all the linear PDE solves in the trace approximation share the same left hand side $\partial_{vu}\bar{r}$ (or its adjoint $\partial_{uv}\bar{r}$), and therefore one can amortize the cost of computing an expensive preconditioner or sparse factorization over several solutions of the linearized PDE. It is worth noting that these different approximations can be exploited in a {multifidelity} framework \cite{PeherstorferWillcoxGunzburger18} where the optimal solution of a cheaper approximation is used as initial guess for a more accurate (and more expensive) approximation. For example, in the numerical results section the optimal solution using the MC-corrected formulation are obtained by first solving the linear approximation in \eqref{eq:Jlin}, then solving the quadratic approximation in \eqref{eq:Jquad} and finally \eqref{eq:JlinMC} or \eqref{eq:JquadMC}.



\begin{table}[!htb]
\caption{The computational cost for the evaluation of the cost functional $\cJ$ and its $z$-gradient $D_z \cJ$ in terms of the number of PDE solves for the five cases of approximation of $\bE[Q]$ and $\text{Var}[Q]$. $N$ is the number of samples/bases for the trace estimator; $p$ is the oversampling factor ($p = 5 \sim 10$) in Algorithm \ref{alg:randomizedEigenSolver}; $M$ is the number of samples used in the sample average approximation or the Monte Carlo correction. 
}\label{table:Cost}
\begin{center}
\begin{tabular}{|c|c|c|c|c|c|}
\hline
Computational cost & MC & Linear & Quadratic & Linear+MC & Quadratic+MC \\
\hline
$\cJ$ (state PDE) & M & 1 & 1 & 1+M & 1+M \\
\hline
$\cJ$ (linear PDE) & 0 & 1 & 1+4N+4p & 1 & 1+4N+4p+2M \\
\hline
$D_z \cJ$ (linear PDE) & M  & 2 & 2+2N & 2+M & 2+2N+3M \\
\hline
\end{tabular}
\end{center}
\end{table}

%
%
%
%
%

\section{Numerical experiments}
\label{sec:numerics}
We consider two examples of PDE-constrained optimal control under uncertainty in which the uncertain parameter is the discretization of an infinite-dimensional random field. These examples aim to demonstrate the scalability, computational efficiency and accuracy of the approximation based on the Taylor approximation and their MC-corrected counterparts. In the first example we consider a fluid flow in a porous medium; the control is the fluid injection rate 
 into the subsurface at specific well locations, the uncertain parameter is the permeability field (modeled as a log-normal random field), and the objective is to drive the pressures measured at the production wells to given target values. In the second example we seek an optimal velocity profile on the inlet boundary of a turbulent jet flow such that the jet width in a given cross-section is maximized; the governing state equations are the Reynolds-averaged Navier--Stokes equations (RANS) with an algebraic closure model and an ad-hoc nonlinear advection-diffusion equation with uncertain coefficients is introduced to model the intermittency of turbulence at the edges of the jet.

The numerical results presented in this section are obtained using \texttt{hIPPYlib} \cite{VillaPetraGhattas2017}, an extensible software framework for large scale PDE-based deterministic and Bayesian inverse problems. \texttt{hIPPYlib} builds on \texttt{FEniCS} \cite{LoggMardalWells2012} (a parallel finite element element library) for the discretization of the PDE and \texttt{PETSc} \cite{BalayAbhyankarAdamsEtAl15} for scalable and efficient linear algebra operations and solvers. Finally, we use the implementation of the L-BFGS algorithm in the \texttt{Python} library \texttt{SciPy}, specifically the \texttt{fmin\_l\_bfgs\_b} routine (BFGS with limited memory and box constraint), to solve the optimization problem.

\subsection{Optimal source control of fluid flow in a porous medium}

We consider an optimal source control problem---originally presented in \cite{AlexanderianPetraStadlerEtAl17}---motivated by flow in a subsurface porous medium, where the state problem is a linear elliptic equation describing single phase fluid flow in a porous medium:
\beq\label{eq:porousmedium}
\begin{array}{rll}
-\nabla \cdot (e^m \nabla u ) &= F z & \text{ in } D,\\
u &= g & \text{ on } \Gamma_D,\\
e^m \nabla u \cdot \boldsymbol{n} & = 0 & \text{ on } \Gamma_N,
\end{array}
\eeq
where $D$ is an open and bounded physical domain $D \subset \bR^d$ ($d= 2, 3$) with Lipschitz boundary $\partial D$, and $\Gamma_D, \Gamma_N \subset \partial D$ ($\Gamma_D \cup \Gamma_N = \partial D$ and $\Gamma_D \cap \Gamma_N = \emptyset$) denote respectively the Dirichlet and Neumann portion of $\partial D$. 
$u$ is the state  representing the pressure. $m$ is an uncertain parameter representing the logarithm of the permeability field, which is assumed to be Gaussian $\cN(\bar{m}, \cC)$  with mean $\bar{m}$ and covariance $\cC$. 
$F$ is a map from the control $z$ to the source term $Fz$ defined as
\beq
Fz = \sum_{i = 1}^{n_c} z_i f_i,
\eeq
where $f_i$, $i = 1, \dots, n_c$, are the mollified Dirac functions located at the $n_c$ injection wells,
and $z = (z_1, \dots, z_{n_c})^\top \in \bR^{N_c}$ is a control that represents the fluid injection rate.
The control $z$ takes admissible value in $\cZ = [z_{\text{min}}, z_{\text{max}}]^{n_c}$, where $z_{\text{min}}$ and $z_{\text{max}}$ denote the minimum and maximum injection rate.

Let $R_g$ denote the lifting of the boundary data $g$ such that $R_g \in H^1(D)$ and $R_g|_{\Gamma_D} = g$. The weak formulation of problem \eqref{eq:porousmedium} reads: Given $m \in \cM = L^2(D)$  and $z \in \cZ$, find $\mathring{u} = u - R_g \in \cU := \{v \in H^1(D): v |_{\Gamma_D} = 0\}$ such that 
\beq\label{eq:porousmediumweak}
r(\mathring{u},v,m,z) = 0, \quad \forall v \in \cV,
\eeq
where $\cV = \cU$, and
\beq
r(\mathring{u},v,m,z) = \int_D e^m \nabla (\mathring{u}+R_g) \cdot \nabla v dx - \int_D F z \; v dx.
\eeq
The control objective measures the difference between the pressure (computed by solving the state problem \eqref{eq:porousmediumweak}) and a target pressure distribution $\bar{u} = (\bar{u}_1, \dots, \bar{u}_{n_p})$ at the production wells $x^1, \dots, x^{n_p} \in D$. Specifically, 
\beq\label{eq:Qporous}
Q(\mathring{u}) = ||\cB \mathring{u} - \bar{u}||^2_2,
\eeq
where $\cB: \cU \to \bR^{n_p}$ is a pointwise observation operator at the locations $x^1, \dots, x^{n_p}$.

\begin{figure}[!htb]
\begin{center}
\includegraphics[scale=0.1]{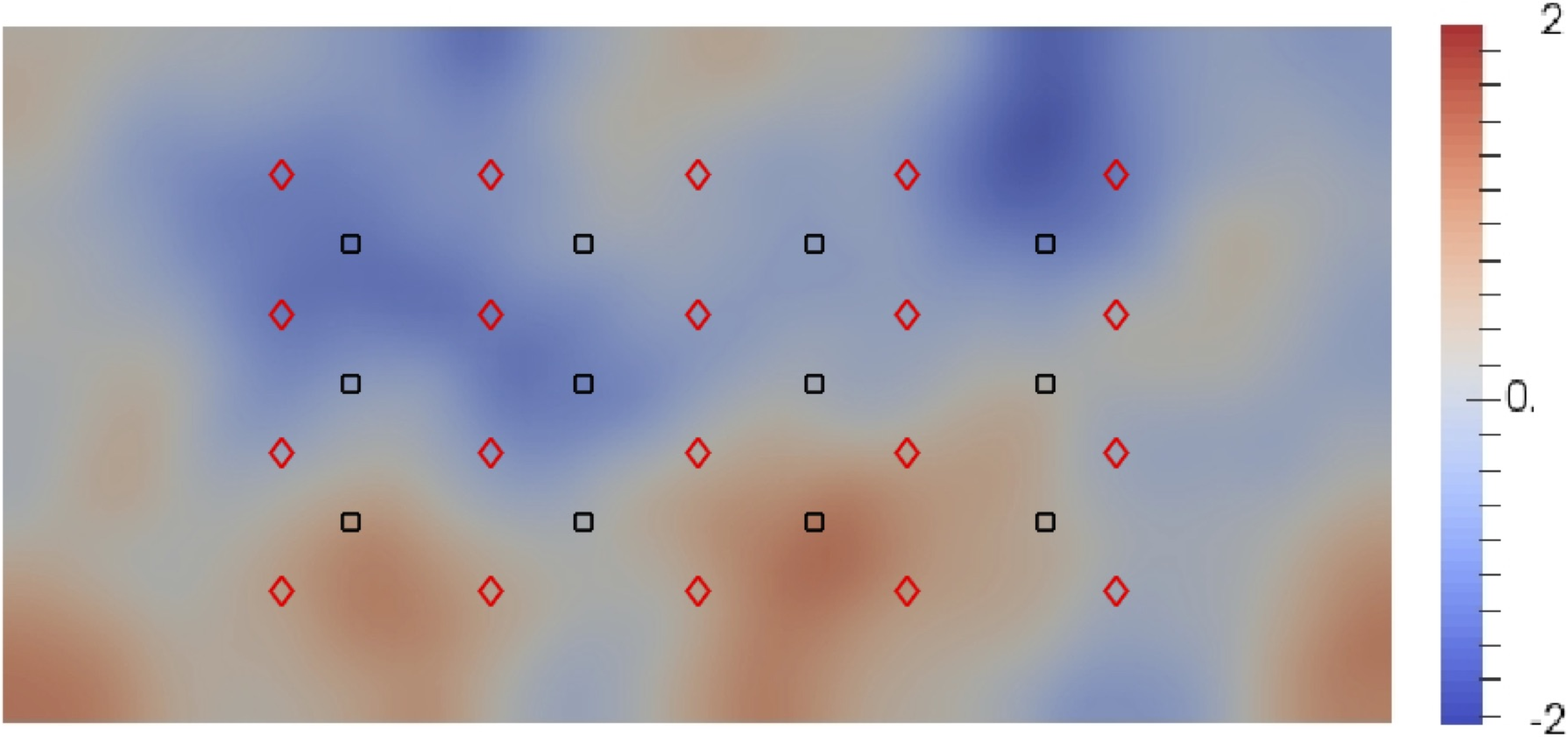}

\includegraphics[scale=0.2]{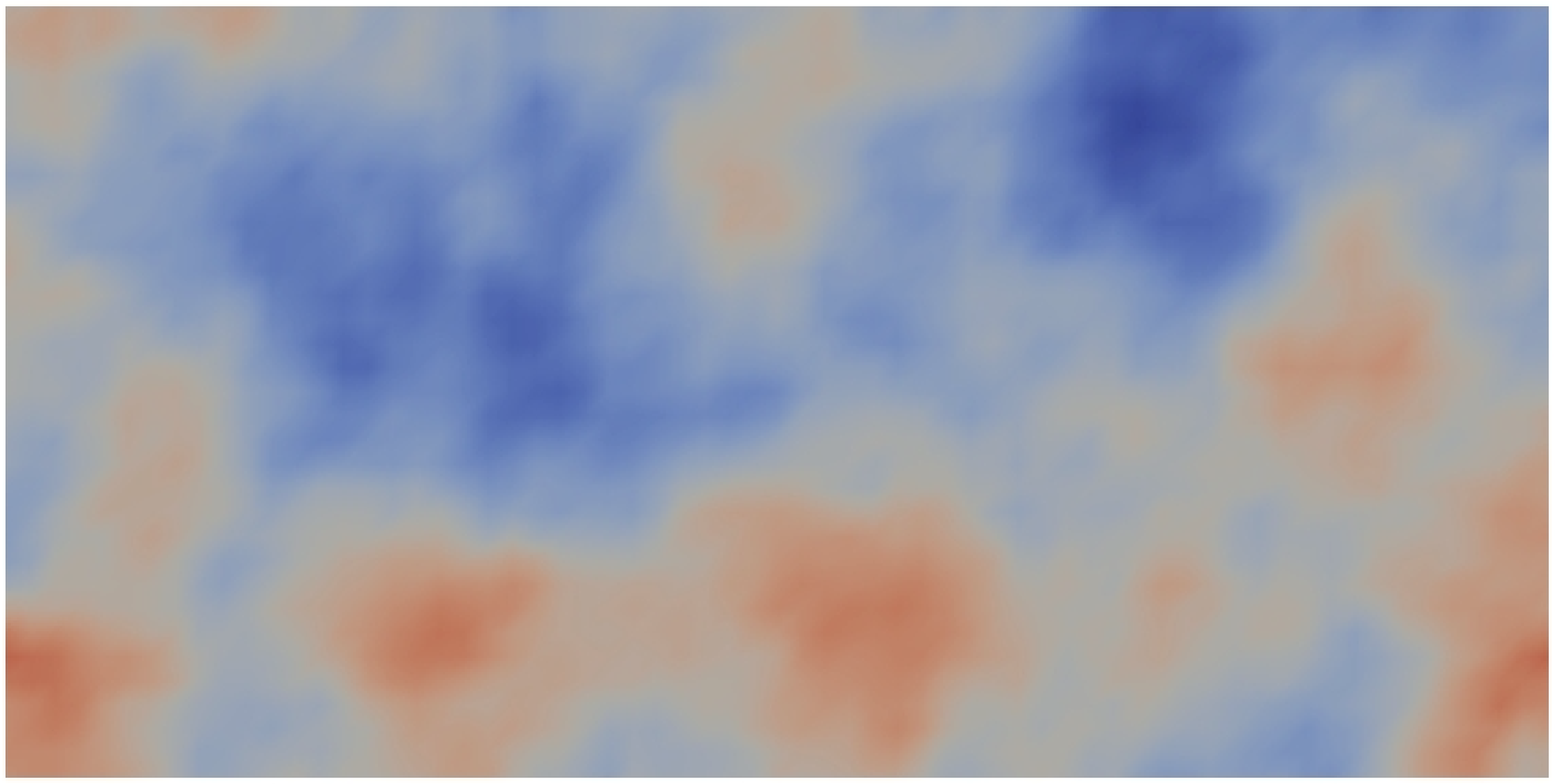}
\includegraphics[scale=0.2]{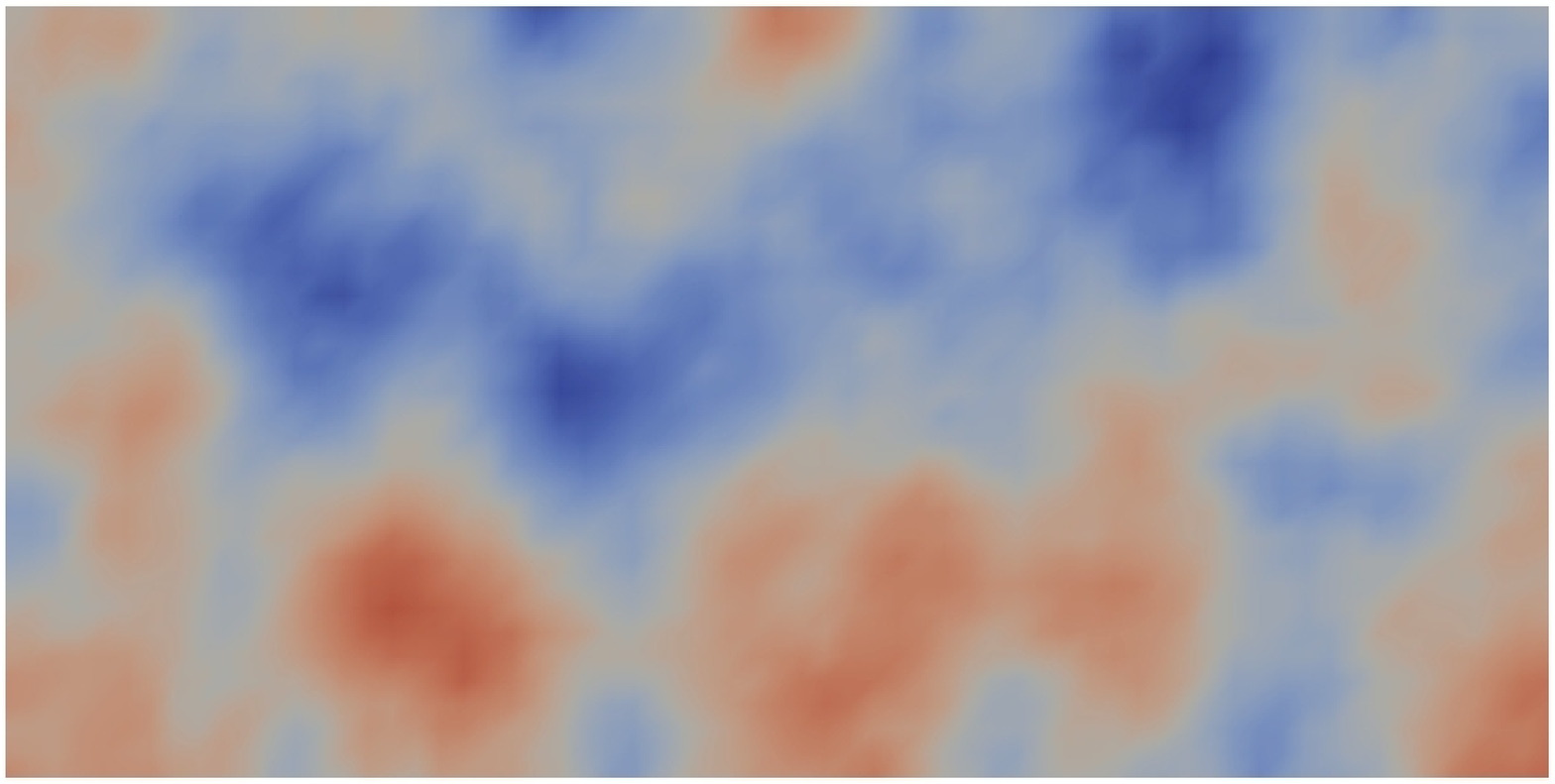}
\includegraphics[scale=0.2]{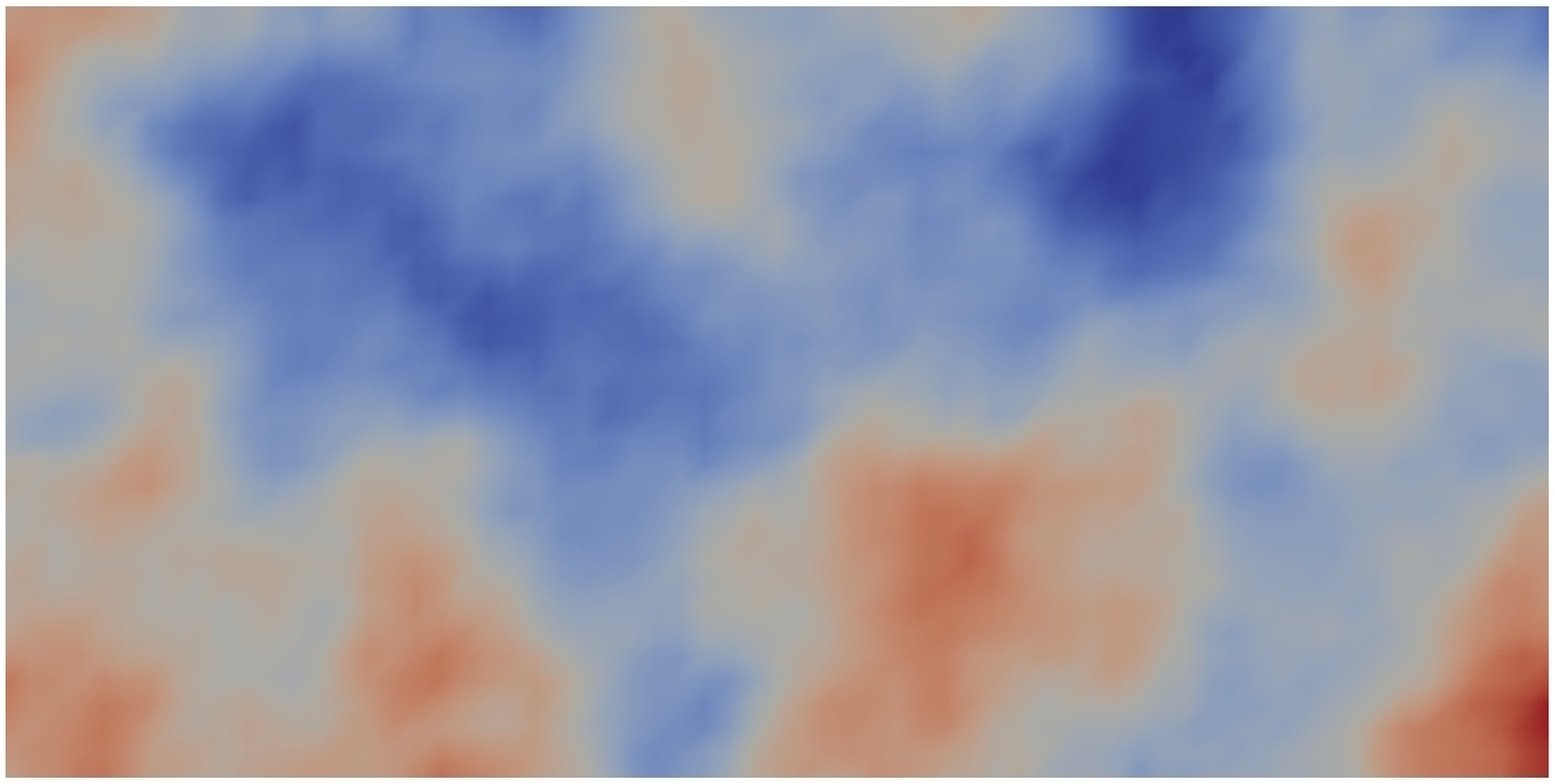}
\end{center}
\caption{Top: mean $\bar{m}$ of the Gaussian  uncertain parameter $m \sim \cN(\bar{m}, \cC)$. The red diamonds represent the locations of the injection wells, and the black squares represent the locations of the production wells. Bottom: random samples drawn from the Gaussian measure $\cN(\bar{m}, \cC)$.}\label{fig:sampleporous}
\end{figure}

Here, we set $D = (0, 2)\times (0, 1)$, $\Gamma_D = \{0, 2\}\times [0, 1]$, and $\Gamma_N = (0, 2)\times \{0, 1\}$. We impose the Dirichlet data $g = 2-x_1$. The mean field $\bar{m}$ is shown in the top part of Fig.\ \ref{fig:sampleporous}, in which locations for the 20 injection wells and the 12 production wells are also shown. For the covariance $\cC$ defined in \eqref{eq:covariance}, we specify $\Theta = [1,0;0,1]$, $\alpha = 2$, $\alpha_1 = 0.1$, and $\alpha_2 = 20$. Three random samples of the Gaussian measure $\cN(\bar{m}, \cC)$ are shown in the bottom part of Fig.\ \ref{fig:sampleporous}.
In the cost functional \eqref{eq:objective}, we set $\beta = 1$ and use the control cost $\cP(z) = 10^{-5} \times ||z||_2^2$. We use a finite element discretization with a structured uniform triangular mesh with $65\times 33$ vertices and piecewise linear elements for the discrete approximation of the state $u$, the adjoint $v$, and the uncertain parameter $m$, which results in 2145 dimensions for $m$. We specify the desired pressure $\bar{u}_i = 3 - 8(x_1^i-1)^2 - 4(x_2^i-0.5)^2$ at the production wells $x^i = (x_1^i, x_2^i)$, $i = 1, \dots, 12$, which are shown in the top part of Fig.\ \ref{fig:sampleporous}. Moreover, we set the lower and upper bounds $z_{\text{min}} = 0$ and $z_{\text{max}} = 32$ for the control, i.e., $z \in [0, 32]^{20}$.

We consider the five different approximations of the cost functional presented in Sec.\ \ref{sec:taylor}, and use the limited memory BFGS method with box constraints to solve the optimization problem, where the $z$-gradient of the cost functional is computed as in Sec.\ \ref{sec:gradientopt}. 
We chose $z_0 = (16, \dots, 16)$ as an initial guess for the control. We used $M = 100$ random samples from $\cN(\bar{m}, \cC)$ to estimate the mean and variance of the control objective in both the sample average approximation presented in Sec.\ \ref{sec:MCinteg} and the variance reduction in Sec.\ \ref{sec:MCVR}. 

\begin{figure}[!htb]
\begin{center}
\includegraphics[scale=0.31]{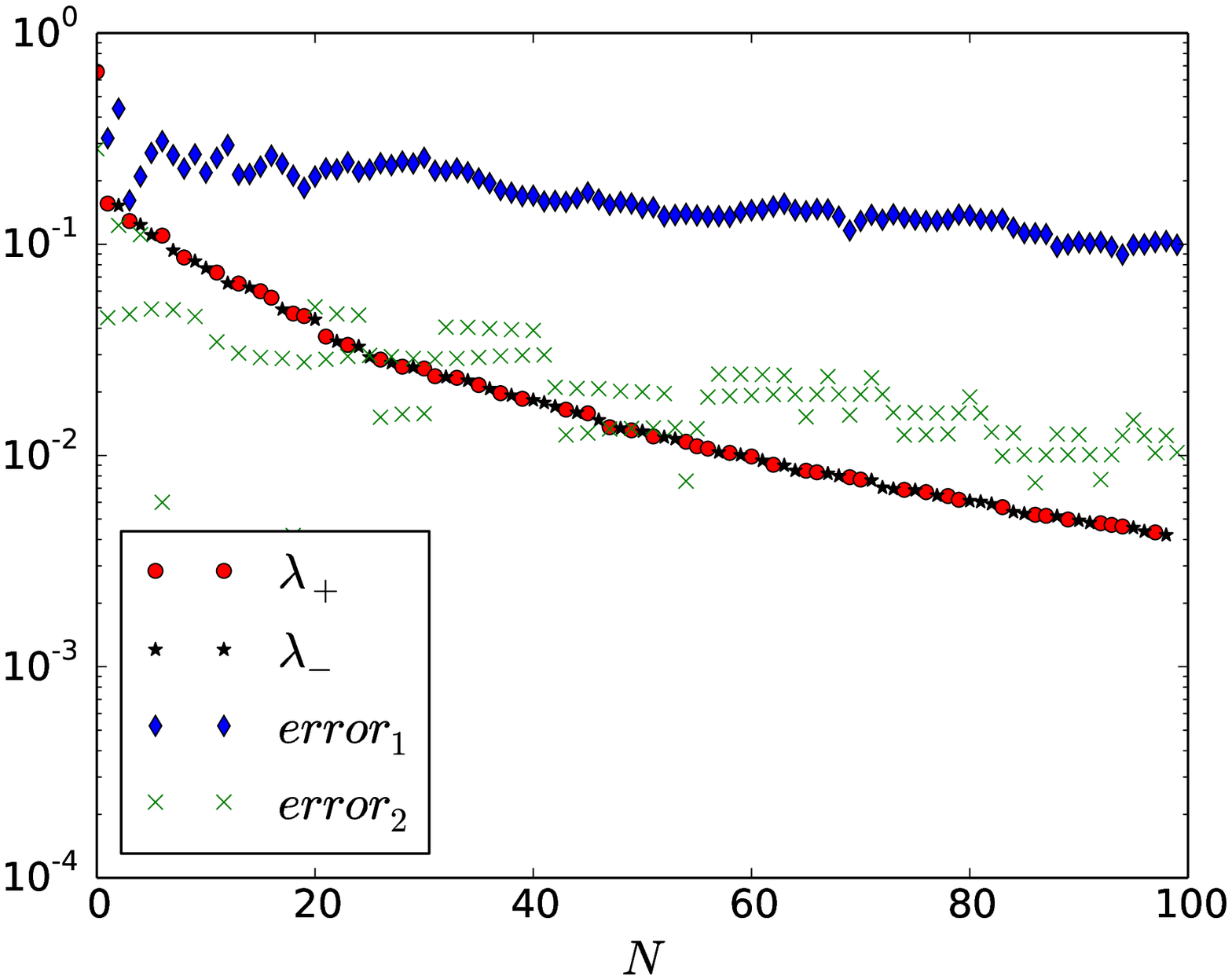}
\includegraphics[scale=0.31]{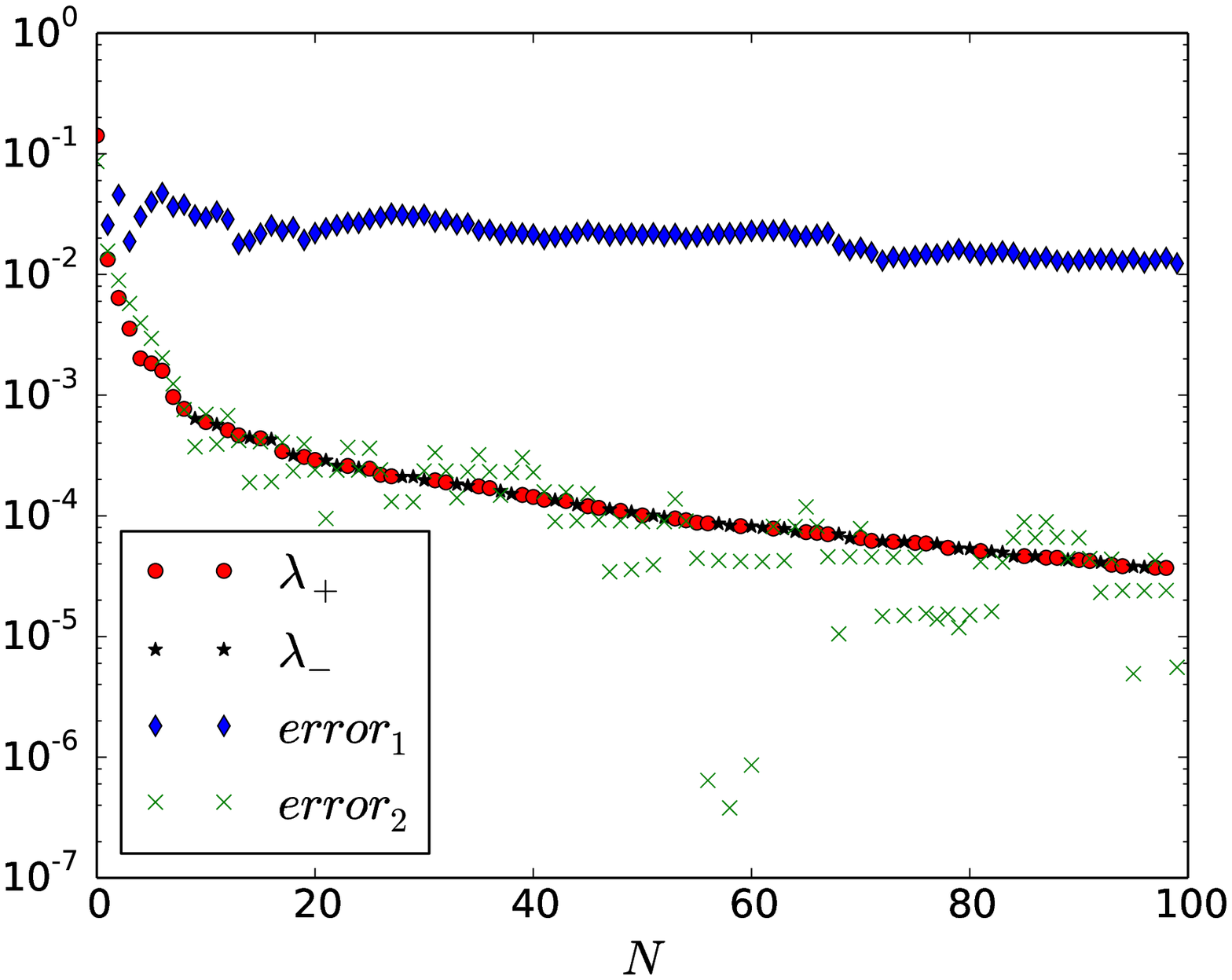}
\end{center}
\caption{The decay of the generalized eigenvalues, $\lambda_+$ for positive eigenvalues and $\lambda_-$ for negative eigenvalues, and the errors, denoted as ${error}_1$ and ${error}_2$ for the trace estimators $\widehat{T}_1$ and $\widehat{T}_2$ with respect to the number $N$ in \eqref{eq:GaussTrace1} and \eqref{eq:RandomizedEigen}. Left: at the initial control $z = z_0$; right: at the optimal control $z = z_{\text{quad}}^{\text{MC}}$ obtained by quadratic approximation with Monte Carlo correction.}\label{fig:traceestimators}
\end{figure}

For the computation of the trace appearing in the quadratic approximations \eqref{eq:quadraticE}, \eqref{eq:quadraticMC}, and \eqref{eq:quadMCVar}, we use the trace estimator $\widehat{T}_2(\cdot)$ \eqref{eq:RandomizedEigen} by solving the generalized eigenvalue problem \eqref{eq:gEigen} using Algorithm \ref{alg:randomizedEigenSolver}, with $k = 100$ and $p = 10$.
We also test the Gaussian trace estimator $\widehat{T}_1(\cdot)$ in \eqref{eq:GaussTrace1}. 
In Fig.\ \ref{fig:traceestimators}, we plot the errors, denoted as $error_1$ and $error_2$ for the two trace approximations, respectively, against the reference value $\sum_{j=1}^{k} \lambda_j$ with $k = 140$ 
at both the initial control $z = z_0$ and the optimal control $z = z_{\text{quad}}^{\text{MC}}$ obtained by the quadratic approximation with Monte Carlo correction. We can observe that in both cases, the error for $\widehat{T}_2(\cH)$ decays much faster than that for $\widehat{T}_1(\cH)$ as a result of fast decaying eigenvalues, especially at the optimal control where more than two orders of magnitude improvement in accuracy is observed for $N = 100$. We remark that the slow decay of the error for $\widehat{T}_1(\cdot)$ may demand a very large $N$ if relatively high accuracy of the trace estimator is required, which leads to a large number of PDE solves at each step of the optimization procedure, thus diminishing the computational advantage of using the Taylor approximation. In the following test, we use $\widehat{T}_2(\cdot)$ with $N = 100$.

\begin{figure}[!htb]
\begin{center}
\includegraphics[scale = 0.31]{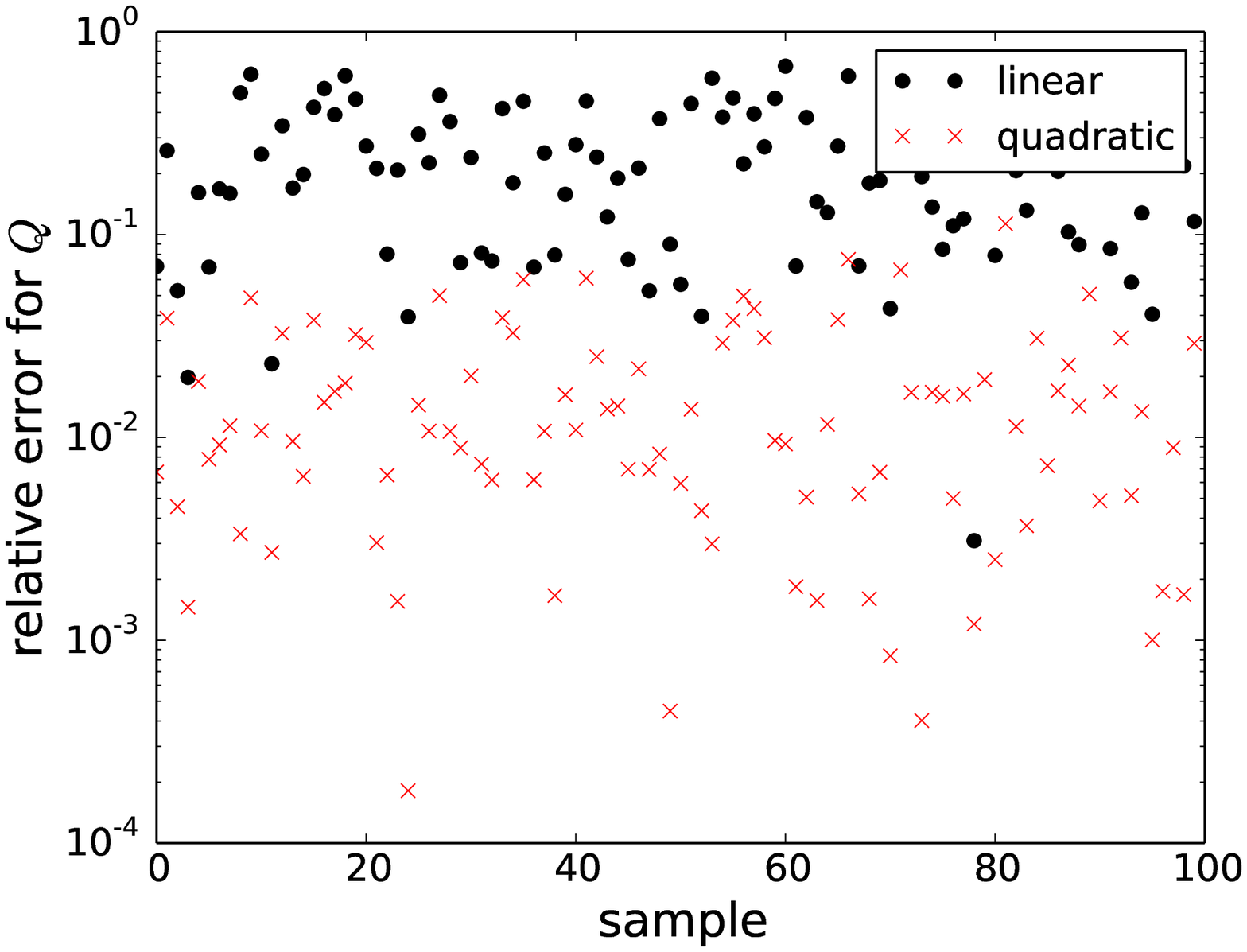}
\includegraphics[scale = 0.31]{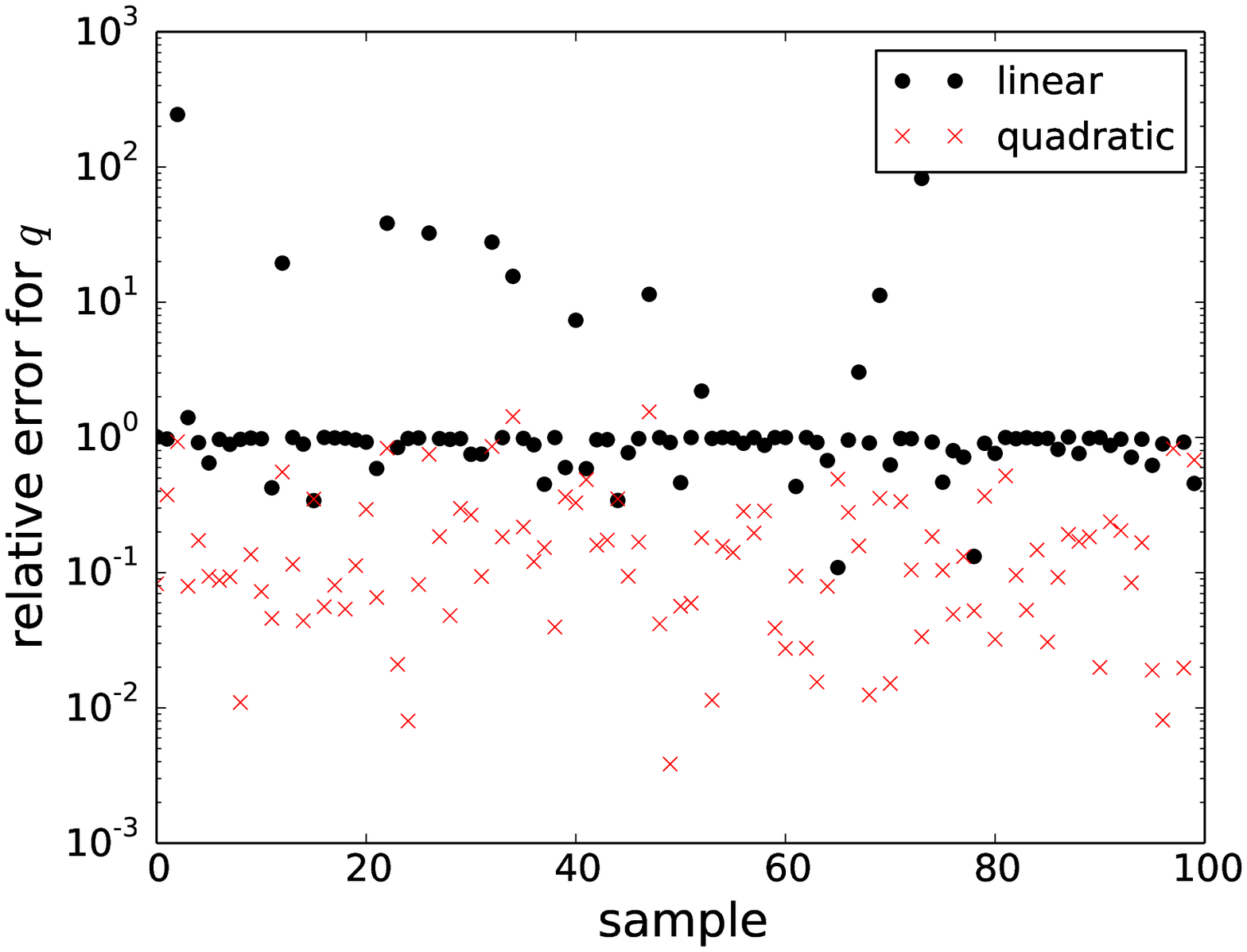}
\end{center}
\caption{Relative errors at 100 random samples of the uncertain parameter for the linear and quadratic approximations of $Q$ used in the computation of the mean $\bE[Q]$ (left) and of $q = (Q-\bar{Q})^2$ for the variance $\text{Var}[Q]$ (right). Here, the errors are show at the optimal control $z = z_{\text{quad}}^{\text{MC}}$ obtained by quadratic approximation with Monte Carlo correction.}\label{fig:linquad}
\end{figure}

Fig.\ \ref{fig:linquad} shows the relative errors at 100 random samples of the uncertain parameter $m$ for the linear and quadratic approximations of $Q$ in \eqref{eq:Qporous} used for the computation of $\bE[Q]$, as well as of $q = (Q-\bar{Q})^2$ in \eqref{eq:linMCVar} and \eqref{eq:quadMCVar} for the computation of $\text{Var}[Q]$. We note that, for both quantities, the quadratic approximation is more accurate than the linear approximation for most of the realizations.

Since the quadratic approximation is very close and strongly correlated to the true value, we may expect large variance reduction by using it as a control variate. In fact, the correlation coefficient between the quadratic approximation and the true value is very close to one, 0.996 for $Q$ and 0.998 for $q$ evaluated using the 100 samples. On the contrary, the linear approximation is poorly correlated to the true value (with correlation coefficient -0.211 for $Q$ and 0.718 for $q$), thus leading to negligible variance reduction.
In Tables \ref{tab:VRmean} and \ref{tab:VRvar}, we report the effect of the variance reduction by the linear and quadratic approximations for the computation of the mean $\bE[Q]$ and variance $\text{Var}[Q]$, respectively. This reduction is measured by the mean-square-error (MSE), which is defined (e.g., for $Q$) as 
\beq\label{eq:MSE}
\text{MSE}(Q) = \frac{\text{Var}_p[Q]}{M},
\eeq
 where $\text{Var}_p[Q]$ is the population variance of $Q$ at $M$ random samples. 
 





\begin{table}[!htb]
\caption{Variance reduction by the linear and quadratic approximations at the initial control $z_0$, and the optimal controls $z_{\text{lin}}^{MC}$ and $z_{\text{quad}}^{MC}$. The mean-square-error (MSE) with different numbers of random samples for the evaluation of the mean $\bE[Q]$ are reported.}\label{tab:VRmean}
\begin{center}
\begin{tabular}{|c|c|c|c|c|c|}
\hline
control &\# samples & $\bE^{\text{MC}}[Q]$& MSE($Q$) & MSE($Q-Q_{\text{lin}}$) & MSE($Q-Q_{\text{quad}}$)  \\
\hline
\multirow{2}{*}{$z_0$} 
&10   & 2.42e+01 & 1.10e+00  & 7.94e-03 & 1.52e-04  \\
&100 & 2.54e+01 & 2.54e-01   & 4.16e-03 & 7.40e-05   \\
\hline
\multirow{2}{*}{$z_{\text{lin}}^{\text{MC}}$} 
&10   & 3.01e-01  & 1.10e-03   & 1.31e-03 & 1.69e-05\\
&100 & 2.89e-01  & 9.68e-05   & 7.49e-05 & 9.00e-07 \\
\hline
\multirow{2}{*}{$z_{\text{quad}}^{\text{MC}}$} 
&10   & 3.16e-01   & 1.11e-03   & 8.02e-04  & 1.05e-05  \\
&100 & 2.98e-01   & 1.10e-04   & 9.98e-05  & 1.12e-06  \\
\hline
\end{tabular}
\end{center}
\end{table}

\begin{table}[!htb]
\caption{Variance reduction by the linear and quadratic approximations at the initial control $z_0$, and the optimal controls $z_{\text{lin}}^{MC}$ and $z_{\text{quad}}^{MC}$. The MSE for the evaluation of the variance $\text{Var}[Q]$, where $q = (Q-\bar{Q})^2$, $q_{\text{lin}} = (Q_{\text{lin}}-\bar{Q})^2$, $q_{\text{quad}} = (Q_{\text{quad}}-\bar{Q})^2$ are reported.}\label{tab:VRvar}
\begin{center}
\begin{tabular}{|c|c|c|c|c|c|}
\hline
control &\# samples & $\bE^{\text{MC}}[q]$& MSE($q$) & MSE($q-q_{\text{lin}}$) & MSE($q-q_{\text{quad}}$)  \\
\hline
\multirow{2}{*}{$z_0$} 
&10   & 1.32e+01 &  1.71e+01 & 9.38e-01  & 6.05e-03  \\
&100 & 2.54e+01 & 1.35e+01  & 2.61e+00 & 3.09e-02 \\
\hline
\multirow{2}{*}{$z_{\text{lin}}^{\text{MC}}$} 
&10   & 1.99e-02  & 1.00e-04   & 9.35e-05  & 3.92e-06 \\
&100 & 1.65e-02  & 1.32e-05   & 1.22e-05  & 2.44e-07 \\
\hline
\multirow{2}{*}{$z_{\text{quad}}^{\text{MC}}$} 
&10   & 2.38e-02  & 1.82e-04   & 1.79e-04  & 8.07e-07  \\
&100 & 1.99e-02  & 1.66e-05   & 1.65e-05  & 3.81e-07  \\
\hline
\end{tabular}
\end{center}
\end{table}

In all cases, we can see that the quadratic approximation achieves more significant variance reduction than the linear approximation. This means that the number of samples needed to achieve a target MSE is significantly smaller for the quadratic approximation than for the linear. 
For example, at the initial control $z_0$, the quadratic approximation results in a speed-up (defined as the ratio between the MSE of the Monte Carlo estimator and the MSE of the correction) of order $1000$X and $100$X in the estimation of $\bE[Q]$ and $\text{Var}[Q]$, respectively. On the other hand, the linear approximation provides a speed-up of order $100$X and $10$X respectively.
In addition, at the optimal controls $z_{\text{lin}}^{\text{MC}}$ and $z_{\text{quad}}^{\text{MC}}$, the quadratic approximation can still achieve a $\sim 100$X speed up for evaluation of both the mean and variance, while the linear approximation is too poor to achieve meaningful variance reduction. Note that the mean and variance are both much smaller at the optimal control than those at the initial control, confirming that the cost functional is significantly reduced by solving the minimization problem \eqref{eq:optimization}.
 
In this example, we see that the trace estimator using the generalized eigenvalues computed by Algorithm \ref{alg:randomizedEigenSolver} is much more accurate than that using the Gaussian randomized (Monte Carlo) trace estimator thanks to the fast decay of the generalized eigenvalues. Moreover, using the Taylor approximation, in particular the quadratic approximation, as a control variate results in over two orders of magnitude in variance reduction for the evaluation of the cost functional.  

\subsection{Optimal boundary control of a turbulent jet flow}
\begin{figure}[!htb]
\begin{multicols}{2}
\vspace*{-0.6cm}
\begin{center}
\includegraphics[width=0.5\textwidth]{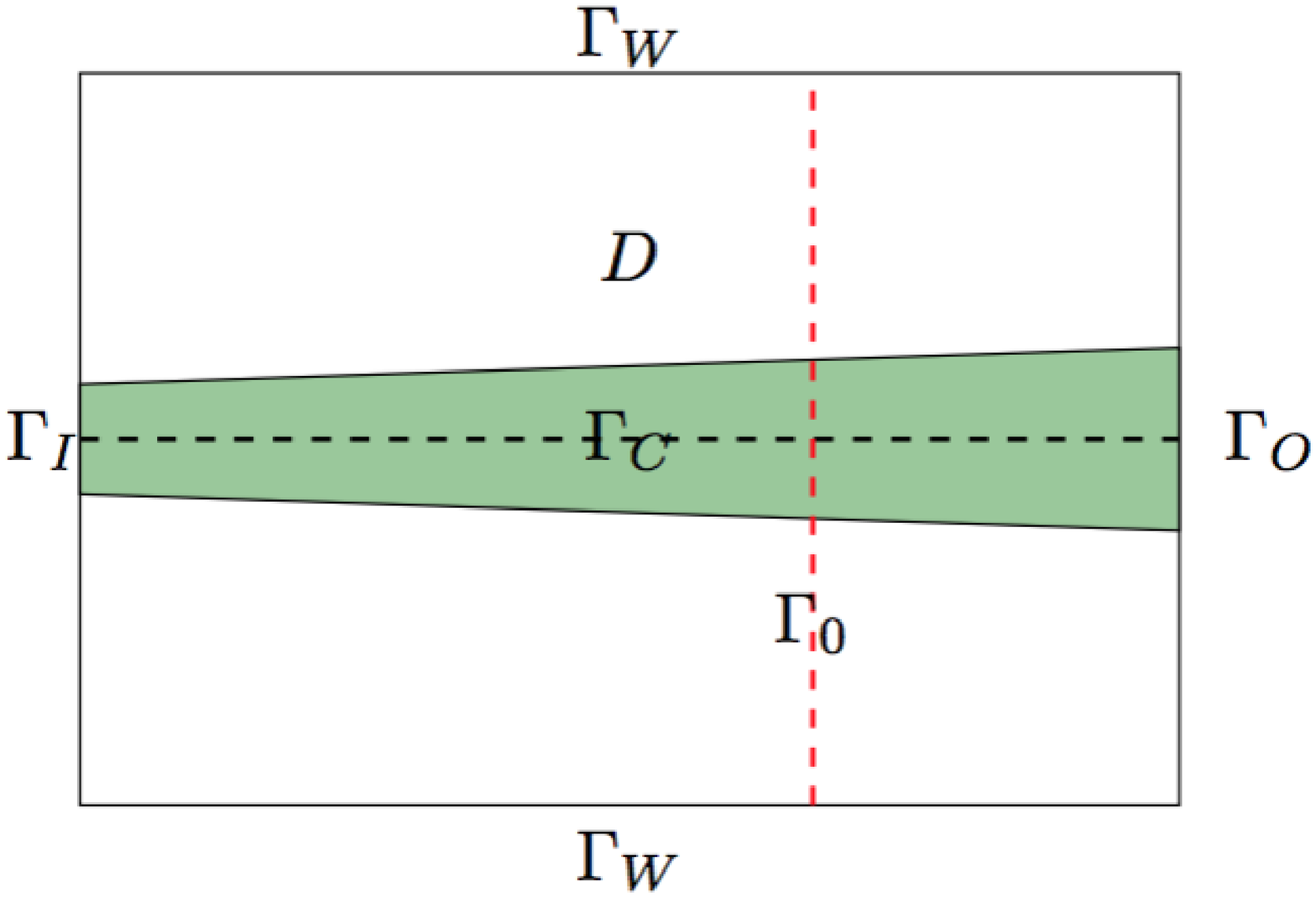}
\end{center}

\columnbreak

\begin{center}
\includegraphics[width=0.45\textwidth]{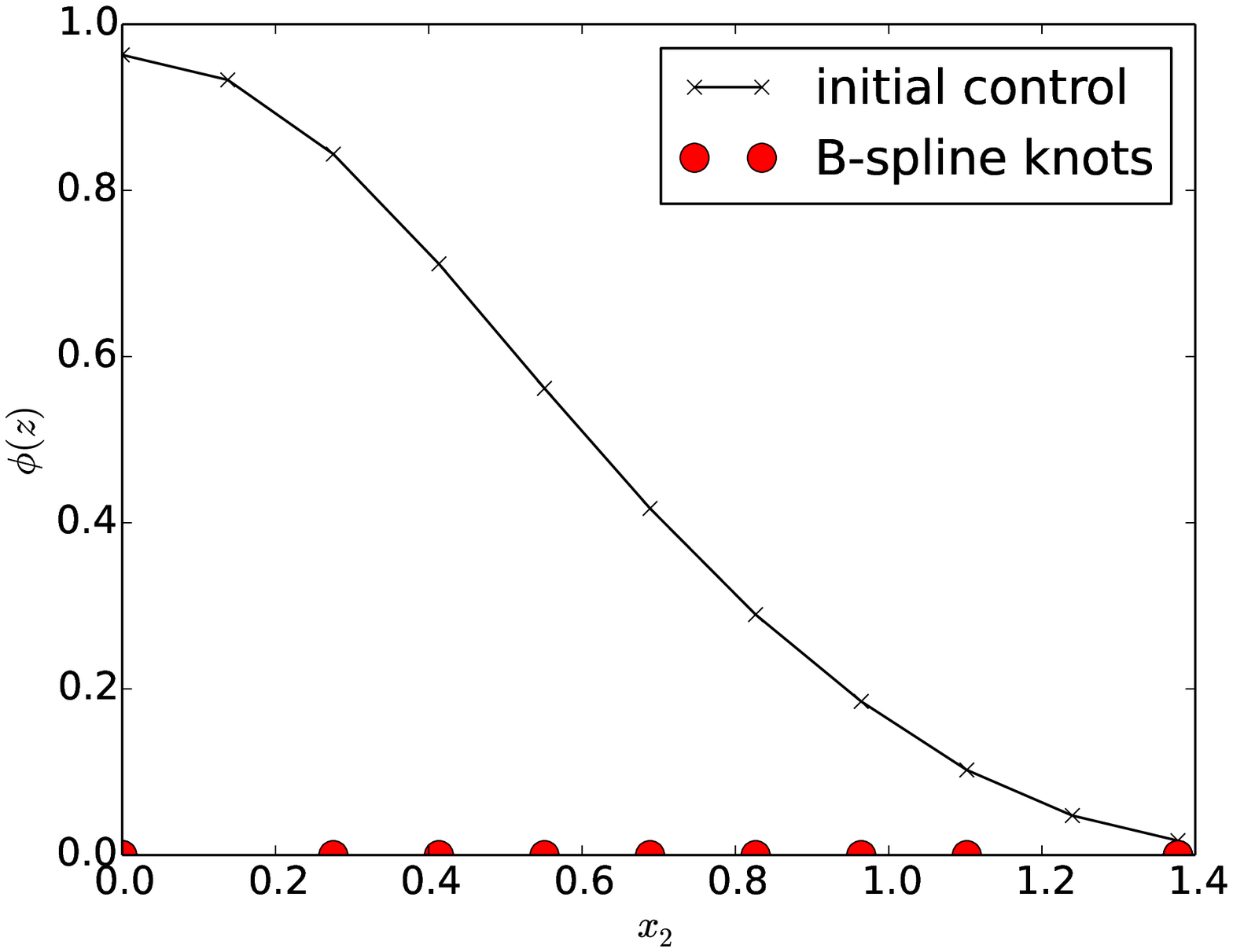}
\end{center}
\end{multicols}
\caption{Left: sketch of the physical domain of the turbulent jet flow, with inlet boundary $\Gamma_I$, outlet boundary $\Gamma_O$, top and bottom wall $\Gamma_W$, center axis $\Gamma_C$, and the cross-section $\Gamma_0$. The computational domain $D$ is the symmetric top half of the physical domain. Right: B-spline interpolation of the inlet velocity profile used as an initial guess, and the distribution of the B-spline knots.}\label{fig:domainsample_t}
\end{figure}

\begin{figure}[!htb]
\begin{center}
\includegraphics[width=0.45\textwidth]{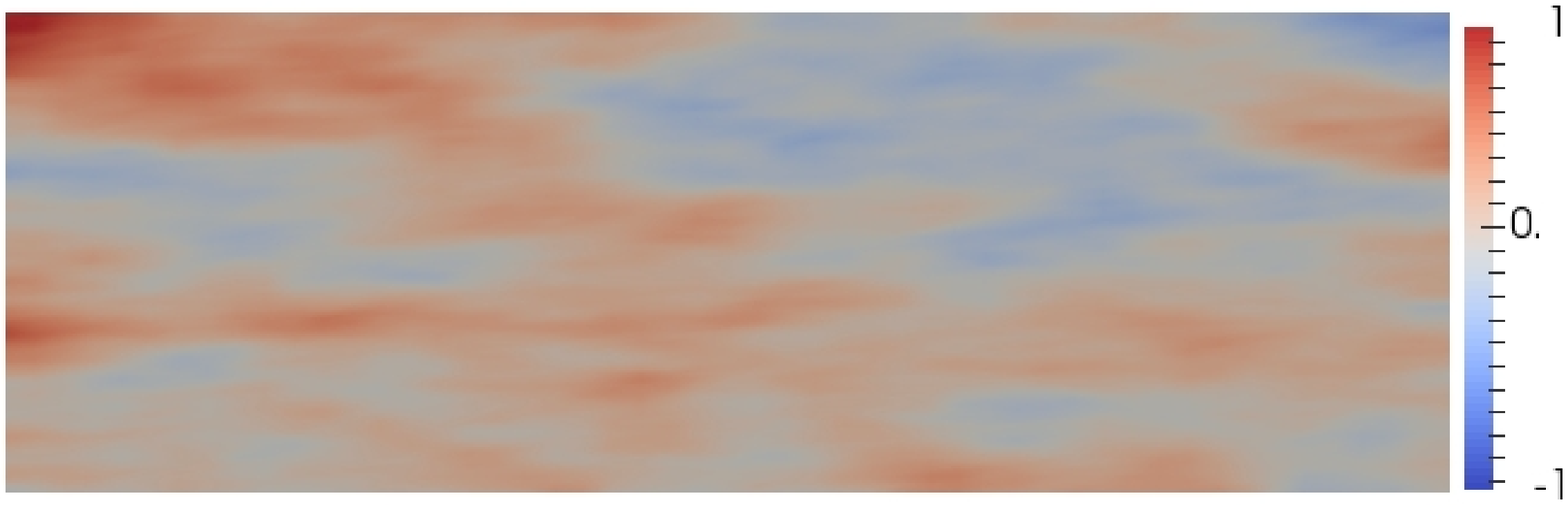}
\includegraphics[width=0.45\textwidth]{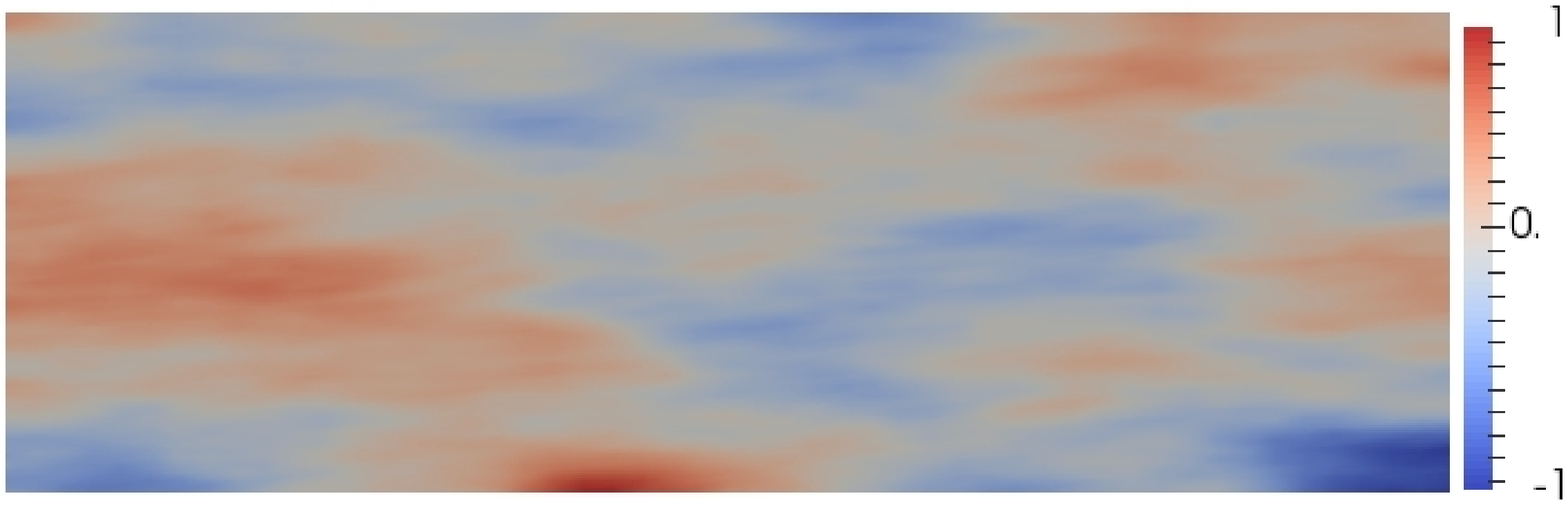}
\end{center}
\caption{Two random samples drawn from the Gaussian distribution $\cN(0, \cC)$ with $\cC$ given in \eqref{eq:covariance}, where $\alpha = 2$, $\alpha_1 = \alpha_2 = 0.5$, $\Theta = [15,0;0,1]$.}\label{fig:sampleN}
\end{figure}

We consider a turbulent jet flow modeled by the Reynolds-averaged Navier--Stokes equations coupled with a nonlinear stochastic advection-diffusion equation that serves as a notional model of uncertainty in the turbulent viscosity $\nu_t$. Specifically, the governing equations read 
\beq\label{eq:turbulenceEq}
\begin{aligned}
- \nabla \cdot \left( \left(\nu +  \nu_t \right) \left( \nabla{\bsu} + \nabla{\bsu}^\top \right)\right) + \left({\bsu} \cdot \nabla\right) {\bsu} + \nabla p &= 0, \quad \text{ in } D,\\
\nabla \cdot {\bsu} &= 0, \quad \text{ in } D, \\
 - \nabla \cdot \left( \left(\nu + \left(\gamma + e^{m} \right)\nu_{t,0} \right) \nabla \gamma \right) + {\bsu} \cdot \nabla \gamma - \frac{1}{2}\frac{\bsu \cdot \bse_1}{x_1 + b} \gamma &= 0, \quad \text{ in } D, 
 \end{aligned}
 \eeq
 with the boundary conditions
 \beq
 \begin{aligned}\label{eq:turbulenceEq_bc}
\bssigma_n(\bsu) \cdot \bstau = 0, \quad \bsu \cdot \bsn + \chi_W \phi(z) & = 0, \quad \text{ on } \Gamma_I, \\
\bssigma_n(\bsu) \cdot \bsn = 0, \quad \bsu \cdot \bstau & = 0, \quad \text{ on } \Gamma_O \cup \Gamma_W, \\
\bssigma_n(\bsu) \cdot \bstau = 0, \quad \bsu \cdot \bsn & = 0, \quad \text{ on } \Gamma_C,\\
\gamma -\gamma_0 & = 0,  \quad \text{ on } \Gamma_I \cup \Gamma_W,\\
\bssigma_n^\gamma(\gamma) \cdot \bsn & = 0,  \quad \text{ on } \Gamma_O \cup \Gamma_C.
\end{aligned}
\eeq
 Here, $D = (0, 30)\times (0, 10)$ is the computational domain (the symmetric upper half of the physical domain for the jet flow) as shown in the left part of Fig.\ \ref{fig:domainsample_t}. $\bse_1 = (1, 0)$, $\bstau$ and $\bsn$ are the unit tangential vector and unit normal vector pointing outside the domain $D$ along the boundary $\partial D$, $\bssigma_n(\bsu) = (\nu + \nu_t) \left( \nabla{\bsu} + \nabla{\bsu}^\top \right) \cdot \bsn$ is the traction vector, and $\sigma_n^\gamma(\gamma) = (\nu + (\gamma + e^m)\nu_{t,0}) \nabla \gamma \cdot \bsn$ is the normal flux of $\gamma$. The control portion of the inflow boundary is denoted by the indicator function $\chi_W$ with $\chi_W(x) = 1$ if $x_2 \in [0, W]$ and $\chi_W(x) = 0$ if $x_2 \in (W, 10)$. Here we take $W = 1.5$. 
The state variables $\bsu = (u_1, u_2)$ and $p$ are the time-averaged velocity and pressure, and $\nu_t = \gamma \nu_{t,0}$ is the turbulent viscosity, where $\nu_{t,0}$ represents an algebraic closure model and $\gamma$ is the indicator function that is close to one near the centerline $\Gamma_C$ and zero outside the jet region. Based on a dimensional analysis of the turbulent planar jet, the turbulent viscosity along the jet centerline is given by
\beq
\nu_{t, 0} = C \sqrt{M}(x_1+aW)^{1/2},
\eeq
where $M = \int_{\Gamma_{I}} ||\bsu||^2 ds$ and $W$ are momentum and jet width at the inlet boundary, $C$ and $a$ are two adimensional parameters that need to be calibrated against either experimental or direct Navier-Stokes simulation (DNS) data. The indicator function 
$\gamma$ satisfies the nonlinear stochastic advection-diffusion equation in \eqref{eq:turbulenceEq} with the boundary data 
\beq
\gamma_0 = 0.5-0.5\tanh\left(5\left(\frac{30-x_1}{30}\right)(x_2-1-0.5x_1)\right).
\eeq
The uncertain parameter $m$ is a Gaussian random field with distribution $\cN(\bar{m}, \cC)$, for which we set the mean $\bar{m}$ as a constant field and the covariance as in \eqref{eq:covariance} with $\alpha = 2$, $\alpha_1 = \alpha_2 = 0.5$, and $\Theta = [15,0;0,1]$. Two samples drawn from $\cN(0, \cC)$ are shown in Fig.\ \ref{fig:sampleN}. In general, the parameters $(C, a, \bar{m}, b)$ are also random variables whose distributions can be inferred from experimental or DNS data. Here, for simplicity, we consider only the uncertainty in $m$ and take deterministic values for $(C, a, \bar{m}, b) = (0.012, 5.29, -0.47, 9.58)$.\footnote{These values are the maximum a posteriori estimates obtained by solving a Bayesian calibration problem where the misfit functional is the $L^2$ distance between the velocity of \eqref{eq:turbulenceEq} and the DNS data provided in \cite{KleinSadikiJanicka2003}.} 

The function $\phi(z)$ in \eqref{eq:turbulenceEq_bc} is the horizontal component of the velocity profile on the inlet boundary $\Gamma_I$, which is represented by B-splines as (we use $t$ for $x_2$ here)
\beq
\phi(z)(t) = \sum_{k=1}^K z_k B_{k,n}(t),
\eeq
where the coefficient vector $z := (z_1, \dots, z_K) \in [0, 1]^K$ is the control and $B_{k,n}(t)$ are B-spline basis functions recursively defined as
\beq
B_{k,n}(t) = \frac{t-t_k}{t_{k+n}-t_k} B_{k, n-1}(t) + \frac{t_{k+n+1}-t}{t_{k+n+1}-t_{k+1}} B_{k+1, n-1}(t),
\eeq
where
\beq
B_{k,0}(t) = \left\{ \begin{array}{ll} 1 & \text{if } t_k \leq t < t_{k+1}; \\0 & \text{otherwise}.\end{array}\right.
\eeq
Here we take $K = 11$ and use the B-spline of order $n = 3$ with $15$ knots ($4$ repetitions at each endpoint), as depicted on the right of Fig.\ \ref{fig:domainsample_t}. A B-spline representation of the inflow velocity profile used for the initial control is also shown. 

The control objective is the jet width $Q$ along the cross-section $\Gamma_0$ at $x_1 = 20$ shown in Fig.\ \ref{fig:domainsample_t}, and is defined as
\beq
Q(\bsu) = \frac{1}{\delta^0_1(\bsu)} \int_{\Gamma_0} \bsu \cdot \bse_1 dx_2, 
\eeq 
where $\delta_1^0(\bsu) = \bsu(x^0) \cdot \bse_1$ with $x^0 = (20, 0)$. 
The objective is to maximize the expected value of the jet width $Q$, i.e., to minimize $\bE[-Q]$ subject to the state problem \eqref{eq:turbulenceEq} and the constraint that the inlet momentum $M$ is prescribed. 
Specifically, we consider the cost functional  
\beq\label{eq:objective_t}
\cJ(z) = \bE[-Q] + \beta \text{Var}[Q] + \cP(z).
\eeq
Here we set $\beta = 10^3$ to balance the mean and variance, and define the penalty term $\cP(z)$ as 
\beq
\cP(z) = \beta_1 \left(\int_{\Gamma_I} (\phi(z))^2 dx - M \right)^2 + \beta_2 \int_{\Gamma_I} |\nabla \phi(z)|^2 dx, 
\eeq
where the first term penalizes violations of the inlet momentum constraint ($\beta_1 = 10^3$), and the second term controls the regularity of the inlet velocity ($\beta_2 = 1$). 

The weak formulation of the turbulence model is given by: Given the realization $m \in \cM = L^2(D)$  and control $z \in \cZ = \mathbb{R}^K$, find $u \in \cU$ such that
\beq\label{eq:weak_t}
r(u, v, m, z) = 0, \quad \forall v \in \cV,
\eeq
where we denote the state with $u = (\mathring{{\bsu}}, p, \mathring{\gamma}) \in \cU$ and adjoint with $v = (\bsv, q, \eta) \in \cV$. The space for the state and adjoint variables are $\cU = \cV = U \times P \times V$, where $U = \{\bsv \in (H^1(D))^2: \bsv\cdot \bsn = 0 \text{ on } \Gamma_C \text{ and } \bsv\cdot \bstau = 0 \text{ on } \Gamma_O \cup \Gamma_W\}$, $P = L^2(D)$, and $V = \{\gamma \in H^1(D): \gamma|_{\Gamma_W \cup \Gamma_I} = 0\}$. Let $R_{0}$ be the lifting of the Dirichlet boundary data for velocity on $\Gamma_C\cup\Gamma_O\cup \Gamma_W$, we have $\bsu = \mathring{\bsu} + R_0$. Let $R_{\gamma_0}$ be the lifting of the boundary data $\gamma_0$ on $\Gamma_I \cup \Gamma_W$, we have $\gamma = \mathring{\gamma}+ R_{\gamma_0}$.

Specifically, the weak form in \eqref{eq:weak_t} consists of three terms and reads
\beq
r(u, v, m, z) = \text{Model}(u,v,m) + \text{Stablization}(u,v) + \text{Nitsche}(u,v,m,z).
\eeq
The first term represents the PDE model in weak form given by
\beq
\begin{aligned}
 \text{Model}(u, v, m) &= \int_D (\nu + \nu_t) 2\bS(\bsu)\cdot \bS(\bsv) dx + \int_D [(\bsu\cdot \nabla) \bsu]\cdot \bsv dx - \int_D p \nabla \cdot \bsv dx + \int_D q \nabla \cdot \bsu dx \\
& + \int_D \left(\nu + \left(\gamma + e^{m} \right)\nu_{t,0} \right) \nabla \gamma \cdot \nabla \eta dx + \int_D [\bsu\cdot \nabla \gamma] \eta dx - \int_D \frac{1}{2} \frac{\bsu\cdot \bse_1}{x_1 + b} \gamma \eta dx.
\end{aligned}
\eeq
where $\bS(\bsu) = (\nabla \bsu + \nabla \bsu^\top)/2$ denotes strain tensor.

The second term of \eqref{eq:weak_t} represents stabilization by the {Galerkin Least-Squares} (GLS) for the momentum and mass conservation equations, and by streamline diffusion for the nonlinear advection diffusion equation, i.e.,
\beq
\begin{aligned}
\text{Stablization}(u,v) &= \int_D \tau_1 L_1(u) \cdot D_u L_1(u)(v) dx  + \int_D \tau_2 (\nabla \cdot \bsu) (\nabla \cdot \bsv) dx  + \int_D \tau_3 (\bsu \cdot \nabla \gamma) (\bsu \cdot \nabla \eta) dx,
\end{aligned}
\eeq
where $L_1(u)$ represents the strong form of the residual of the momentum equation (line 1) of \eqref{eq:turbulenceEq}, and $\tau_1$, $\tau_2$ and $\tau_3$ are properly chosen stabilization parameters associated with the local P\'eclet number.

The third term of \eqref{eq:weak_t} represents the weak imposition of the inlet velocity profile (i.e., the control) by Nitsche's method \cite{BazilevsHughes2007} and reads  
\beq
\begin{aligned}
\text{Nitsche}(u,v,m,z) & =   C_d \int_{\Gamma_I} h^{-1} (\nu + \nu_t) (\bsu \cdot \bsn + \chi_W \phi(z)) (\bsv \cdot \bsn) ds \\
& - \int_{\Gamma_I} (\sigma_n(\bsu) \cdot \bsn) (\bsv \cdot \bsn) +  (\sigma_n(\bsv) \cdot \bsn) (\bsu \cdot \bsn + \chi_W \phi(z)) ds, 
\end{aligned}
\eeq
where the first term guarantees the coercivity of the form and the second term arises from integration by parts.
The constant $C_d$ is taken as $C_d = 10^5$, and $h$ is a local element edge length along the Dirichlet boundaries. This weak imposition of the inlet velocity is necessary to reveal the boundary control in the optimization formulation.

We discretize the problem using finite element with piecewise polynomials (P2, P1, P1) for the state $u$ and P1 for the uncertain parameter $m$. We triangulate the computational domain using an anisotropic mesh, locally refined in the jet region, with 14,400 triangular elements ($60$ elements along the $x_1$-direction and $120$ along the $x_2$-direction). After discretization, we obtain an optimization problem in which the dimension of the control, state, and uncertain parameter are 11, 73084, and 7381, respectively. As in the previous example, we use a BFGS method with bound constraints to solve the optimization problem \eqref{eq:optimization} with the cost functional defined in \eqref{eq:objective_t} and the constraint defined by the state problem \eqref{eq:weak_t}. At each BFGS iteration, we solve the fully-coupled nonlinear state problem using Newton's method with an LU factorization of the Jacobian operator. 
To compute all the derivatives involved in the $z$-gradient presented in Sec.\ \ref{sec:gradientopt} we exploit the symbolic differentiation capabilities of \texttt{FEniCS} \cite{LoggMardalWells2012}.

\begin{figure}[!htb]
\begin{center}
\includegraphics[scale=0.31]{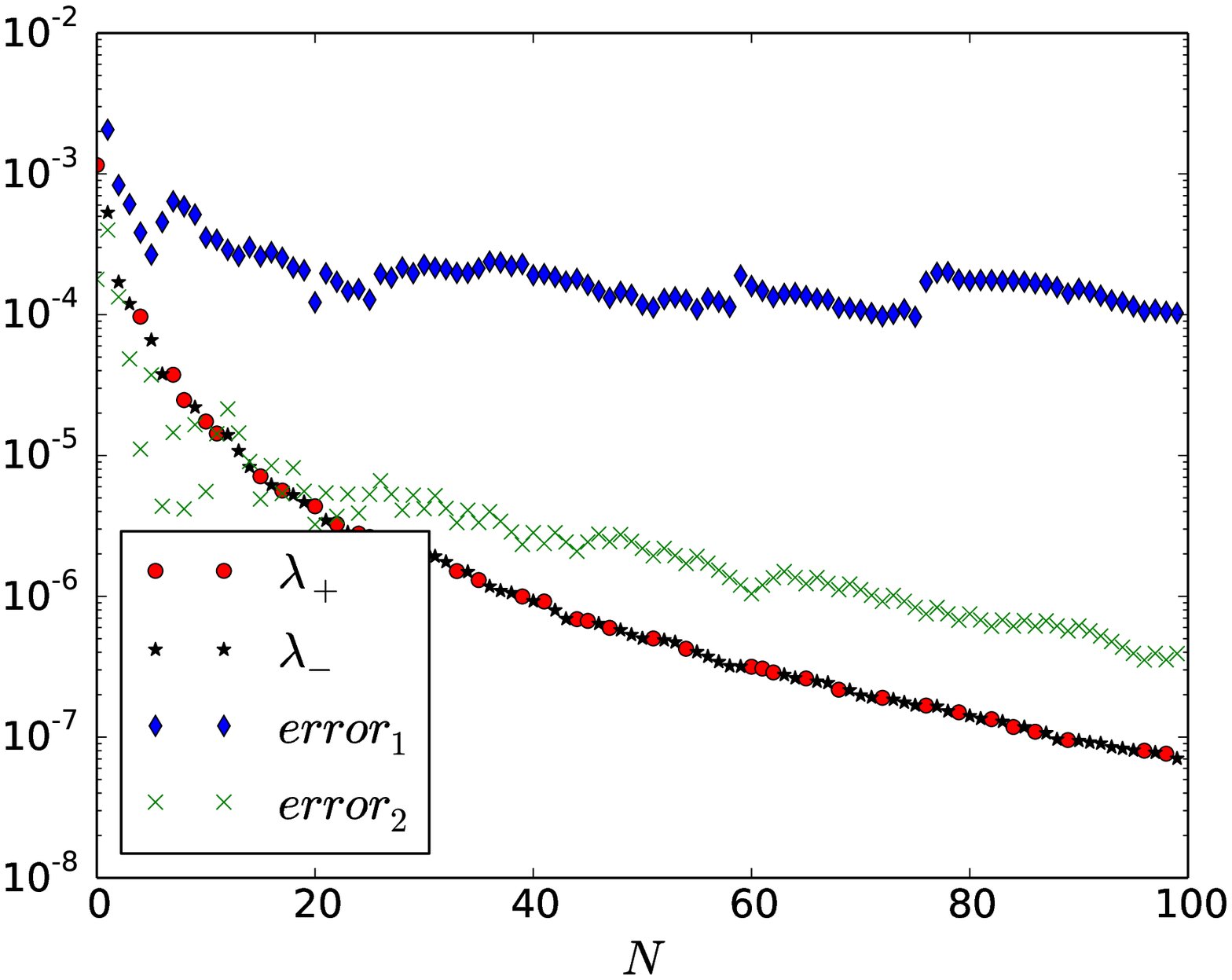}
\includegraphics[scale=0.31]{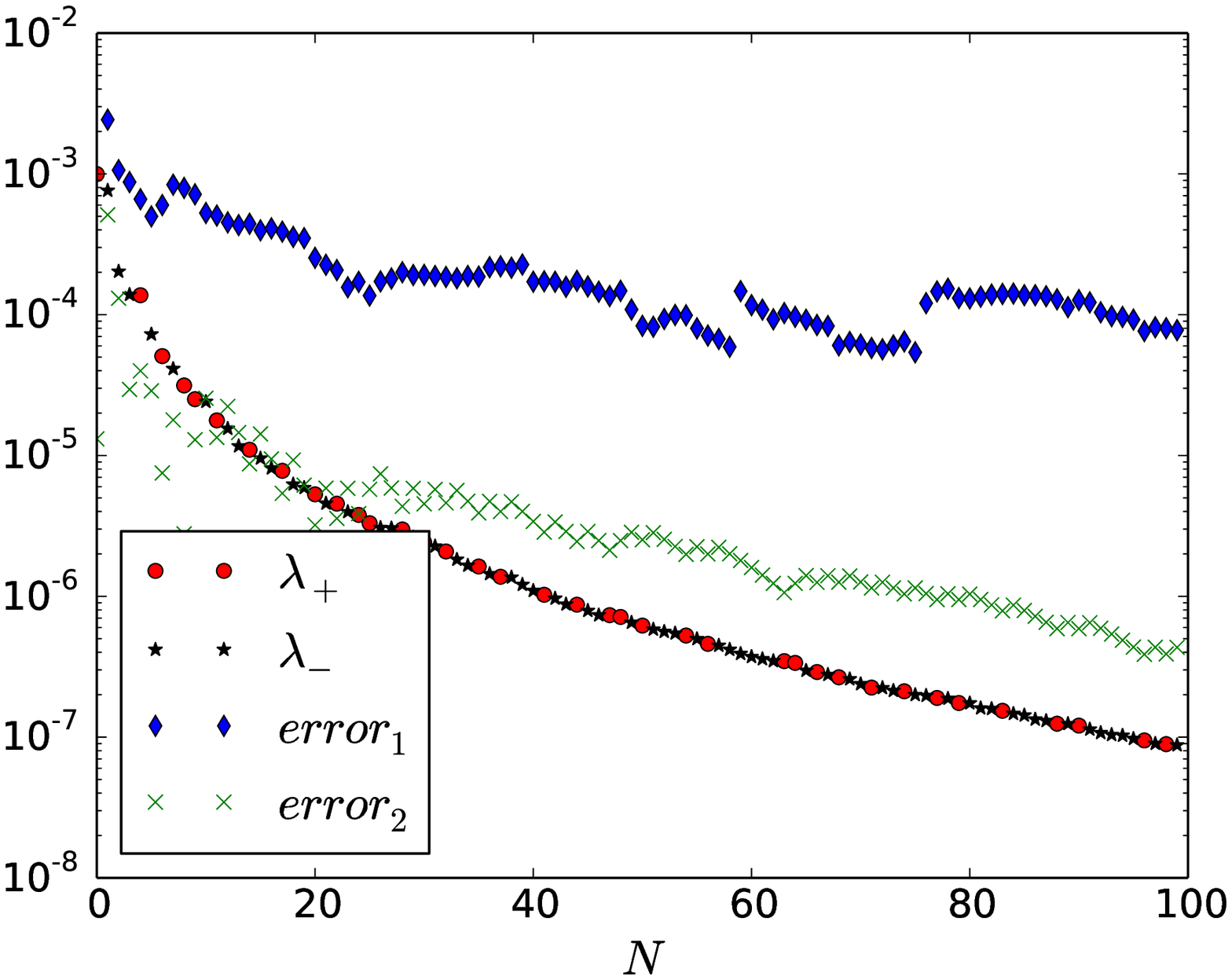}
\end{center}
\caption{The decay of the generalized eigenvalues, $\lambda_+$ for positive eigenvalues and $\lambda_-$ for negative eigenvalues, and the errors, denoted as ${error}_1$ and ${error}_2$ for the trace estimators $\widehat{T}_1$ and $\widehat{T}_2$ with respect to the number $N$ in \eqref{eq:GaussTrace1} and \eqref{eq:RandomizedEigen}. Left: the initial control $z = z_0$; right: the optimal control $z = z_{\text{quad}}^{\text{MC}}$ obtained by quadratic approximation with Monte Carlo correction.}\label{fig:traceestimators_turbulence}
\end{figure}

For the computation of the trace in the quadratic approximation we use both the Gaussian trace estimator $\widehat{T}_1(\cdot)$ in \eqref{eq:GaussTrace1},
and the trace estimator $\widehat{T}_2(\cdot)$ in \eqref{eq:RandomizedEigen} obtained by solving the generalized eigenvalue problem \eqref{eq:gEigen}. Fig.\  \ref{fig:traceestimators_turbulence} displays the decay of the eigenvalues as well as the approximation error of the two estimators for both the initial control $z_0$ and the optimal control $z_{\text{quad}}^{\text{MC}}$. The rapid decay of the eigenvalues (four orders of magnitude reduction in the first 100 out of 7381 eigenvalues) of the covariance-preconditioned Hessian reflects the intrinsic low-dimensionality of the map from the uncertain parameter field $m$ to the control objective $Q$. This is despite the complexity and nonlinearity of the state equations.
The reference value for the trace $\text{tr}(\cH)$ was computed (with a precision of $10^{-8}$) as the sum of the first 140 dominant eigenvalues of $\cH$ using Algorithm \ref{alg:randomizedEigenSolver} with $k = 140$ and $p = 10$. 
It is evident that the second trace estimator achieves much faster decay of the error than the Gaussian trace estimator, with a gain in accuracy of about three orders of magnitude when using $N = 100$ samples/eigenvalues. In the numerical test, we use the second trace estimator with $N = 20$, which already provides sufficiently accurate results. It is worth noticing that, in contrast, the Gaussian trace estimator would have required about $N=10^4$ samples to achieve the same accuracy.

\begin{figure}[!htb]
\begin{center}
\includegraphics[scale=0.4]{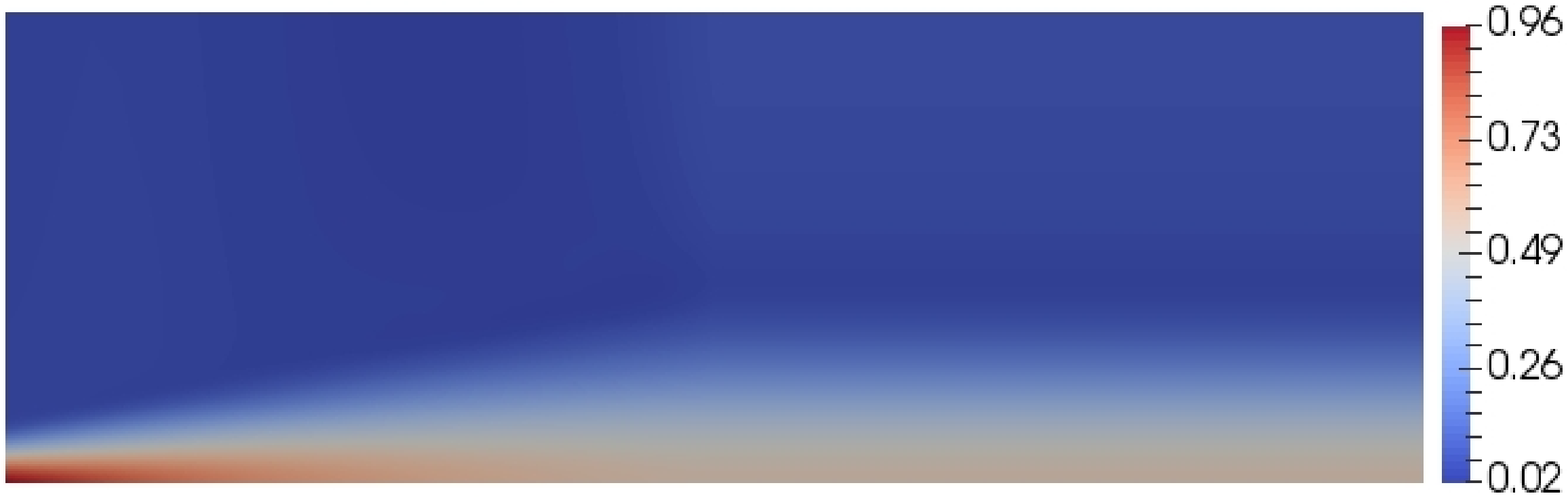}
\includegraphics[scale=0.4]{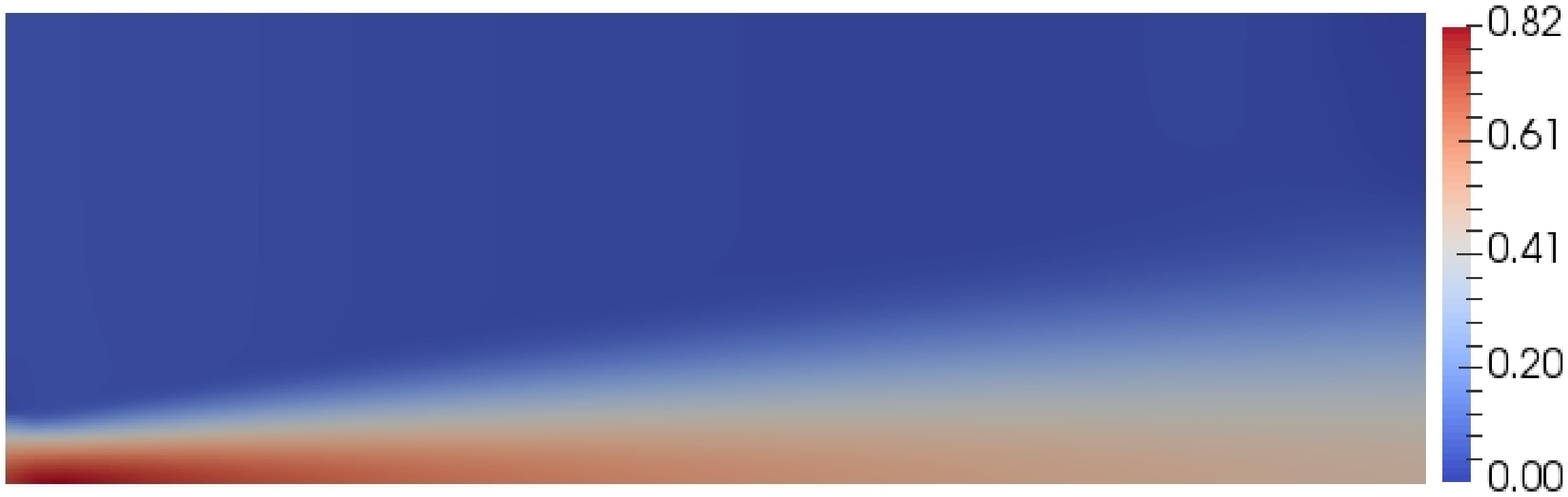}
\end{center}
\caption{The velocity field of the turbulent jet flow obtained solving the state model for the initial (left) and optimal (right) inlet velocity profile. Specifically, the optimal profile was obtained using the quadratic approximation with Monte Carlo correction.}\label{fig:velocitydnscontrol}
\end{figure}

\begin{figure}[!htb]
\begin{center}
\includegraphics[scale = 0.31]{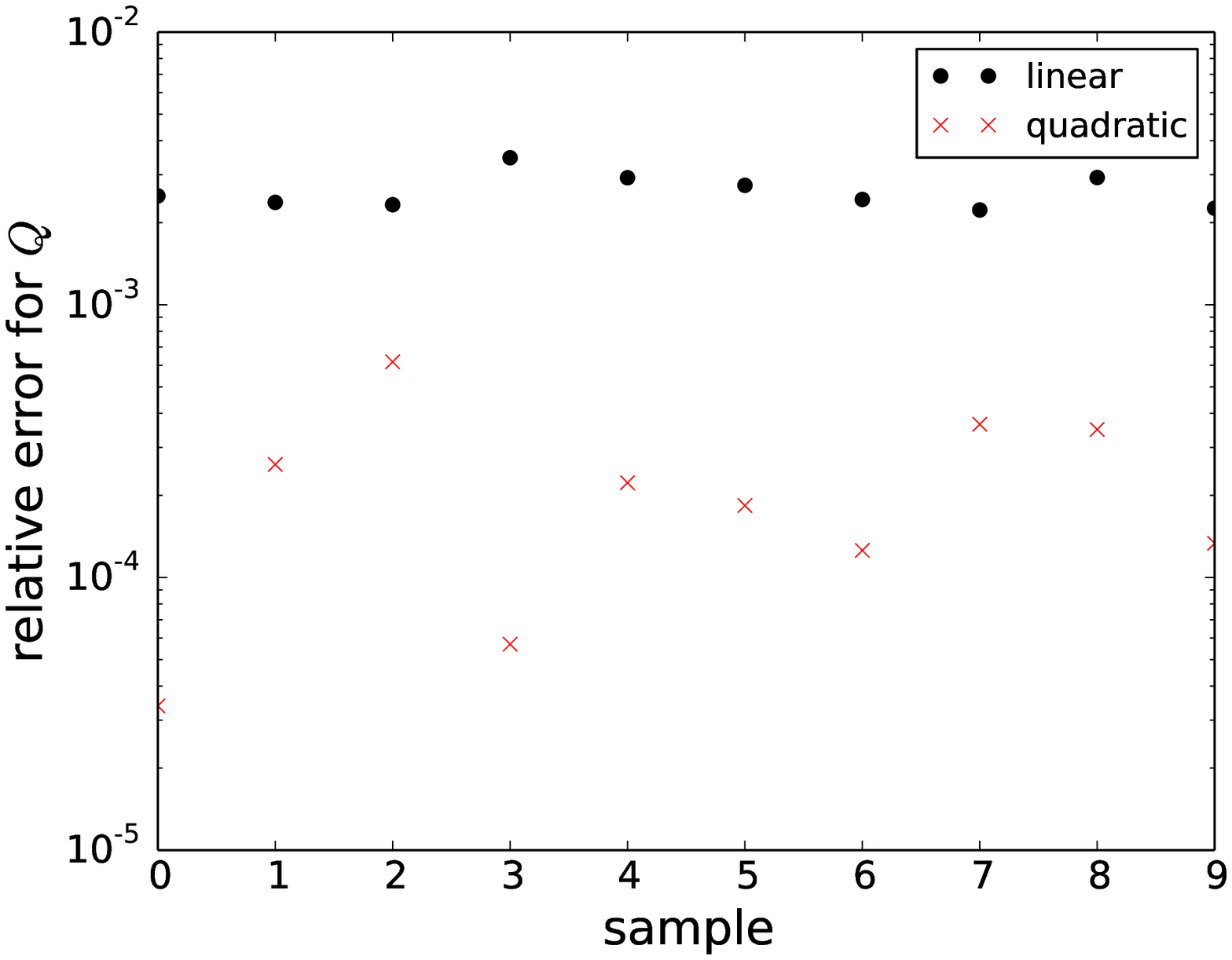}
\includegraphics[scale = 0.31]{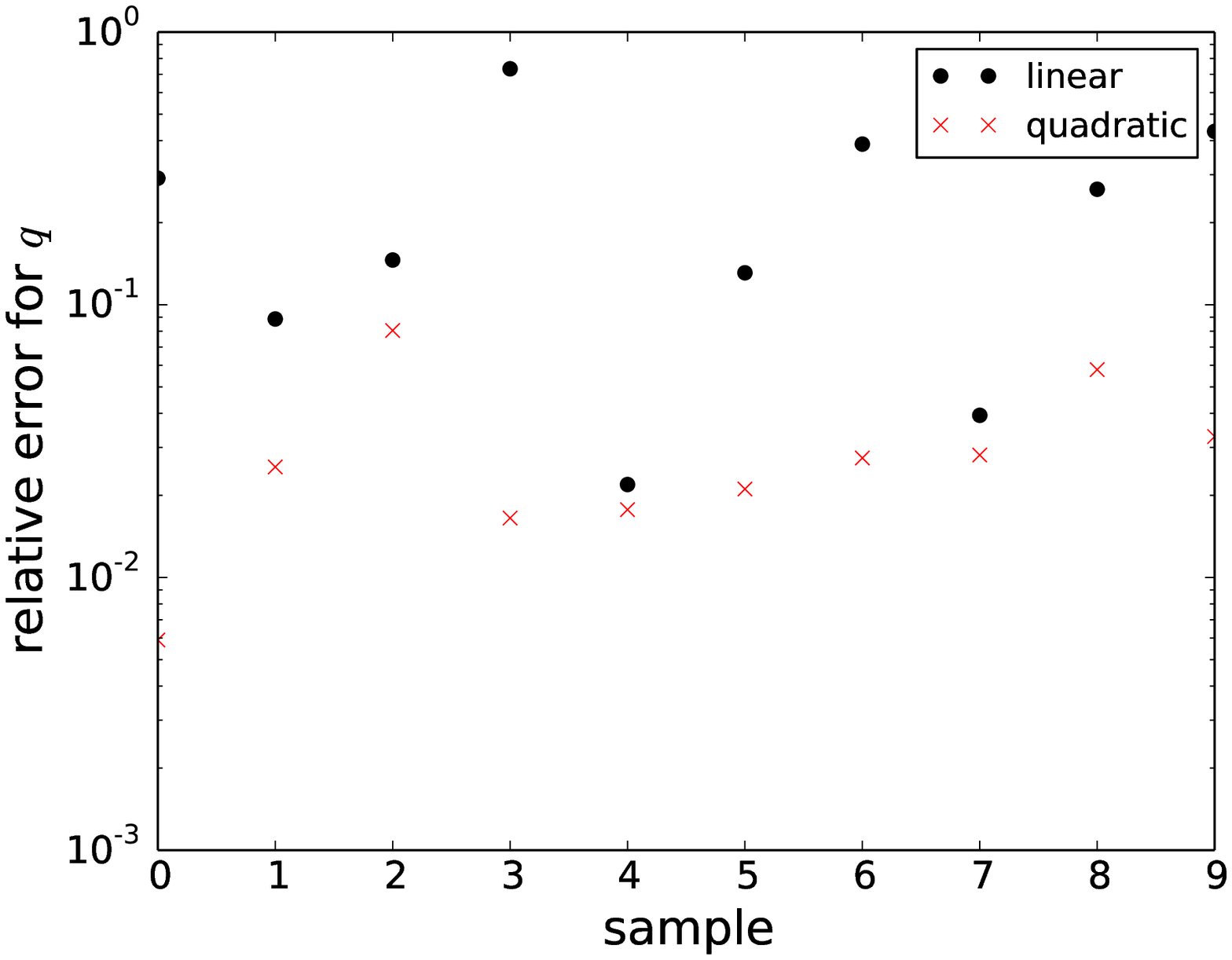}
\end{center}
\caption{Relative errors at 10 random samples of the linear and quadratic approximations of $Q$ (for the computation of the mean $\bE[Q]$ (left)) and of $q = (Q-\bar{Q})^2$ (for the variance $\text{Var}[Q]$ (right)). Here, the errors are shown at the optimal control $z = z_{\text{quad}}^{\text{MC}}$.}\label{fig:linquad_turbulence}
\end{figure}

Fig.\  \ref{fig:linquad_turbulence} shows the relative difference between the true quantity of interest and the linear and quadratic approximations evaluated at 10 random samples for both $Q$ and $q = (Q-\bar{Q})^2$. It is clear that the quadratic approximation is more accurate than the linear approximation for both quantities. Nevertheless both approximations lead to drastic reduction in the variance of the Monte Carlo correction with respect to plain Monte Carlo estimator, as demonstrated by the MSE reported in Tables \ref{tab:VRmean_turbulence} and \ref{tab:VRvar_turbulence}. This variance reduction translates to greater than $100$X speed up in the computation relative to the plain Monte Carlo estimator for the same accuracy (MSE). 

\begin{table}[!htb]
\caption{Variance reduction by the linear and quadratic approximations at the optimal control $z_{\text{lin}}$ for different uncertain parameter dimensions. The mean-square-error (MSE) with 10 random samples for the evaluation of the mean $\bE[Q]$ are reported.}\label{tab:VRmean_turbulence}
\begin{center}
\begin{tabular}{|c|c|c|c|c|c|}
\hline
parameter dimension & $\bE^{\text{MC}}[-Q]$& MSE($Q$) & MSE($Q-Q_{\text{lin}}$) & MSE($Q-Q_{\text{quad}}$)  \\
\hline
1,891  & -1.71e+00 & 7.40e-06 & 2.68e-08 & 1.81e-09   \\
\hline
7,381  & -1.59e+00 & 7.94e-06 & 1.57e-07 & 1.46e-08 \\
\hline
29,161& -1.44e+00 & 3.82e-06 & 7.23e-08 & 1.66e-08   \\
\hline
115,921& -1.42e+00 & 9.47e-06 & 6.91e-08 & 3.06e-08    \\
\hline
462, 241 & -1.41e+00 & 8.85e-06 & 9.00e-08 & 1.78e-08 \\
\hline
\end{tabular}
\end{center}
\end{table}

\begin{table}[!htb]
\caption{Variance reduction by the linear and quadratic approximations at the optimal control $z_{\text{lin}}$ for different uncertain parameter dimensions. The MSE with test 10 random samples for the evaluation of the variance are reported, where $q = (Q-\bar{Q})^2$, $q_{\text{lin}} = (Q_{\text{lin}}-\bar{Q})^2$, $q_{\text{quad}} = (Q_{\text{quad}}-\bar{Q})^2$.}\label{tab:VRvar_turbulence}
\begin{center}
\begin{tabular}{|c|c|c|c|c|c|}
\hline
parameter dimension & $\bE^{\text{MC}}[q]$& MSE($q$) & MSE($q-q_{\text{lin}}$) & MSE($q-q_{\text{quad}}$)  \\
\hline
1,891   & 8.05e-05 & 9.37e-10 & 1.76e-11 & 8.77e-13  \\
\hline
7,381   & 8.13e-05 & 1.15e-09  & 8.87e-12 & 1.48e-12 \\
\hline
29,161  & 5.60e-05 & 6.59e-10 & 3.39e-11 & 4.04e-12   \\
\hline
115,921  & 1.07e-04 & 1.82e-09 & 7.04e-11 & 2.08e-11    \\
\hline
462,241  & 8.96e-05 & 1.37e-09 & 7.78e-12  & 1.33e-11    \\
\hline
\end{tabular}
\end{center}
\end{table}

To investigate the scalability of the proposed algorithm, we consider up to six levels of uniform mesh refinement: $30\times 60$, $60\times 120$,  $120\times 240$,  $240\times 480$,  $480\times 960$, and  $720\times 1440$. The dimension of the uncertain parameter $m$ ranges from $1,891$ on the coarsest mesh to $1,038,961$ on the finest. For the linear and quadratic approximation, we solve the optimization problem using the BFGS algorithm with convergence defined by an absolute tolerance of $10^{-3}$ for the $\ell^{\infty}$ norm of the projected gradient. We also use a mesh continuation technique to define initial guesses for the control: on the coarsest mesh we use the velocity profile in Fig.\  \ref{fig:domainsample_t} (right), and, on the finer meshes, we use the optimal control obtained from the previous (coarser) mesh. For the linear approximation with Monte Carlo correction, we set the BFGS absolute tolerance to $10^{-4}$ and we use the solution of the linear approximation (at the same mesh resolution) as the initial guess.

\begin{figure}[!htb]
\begin{center}
\includegraphics[scale=0.31]{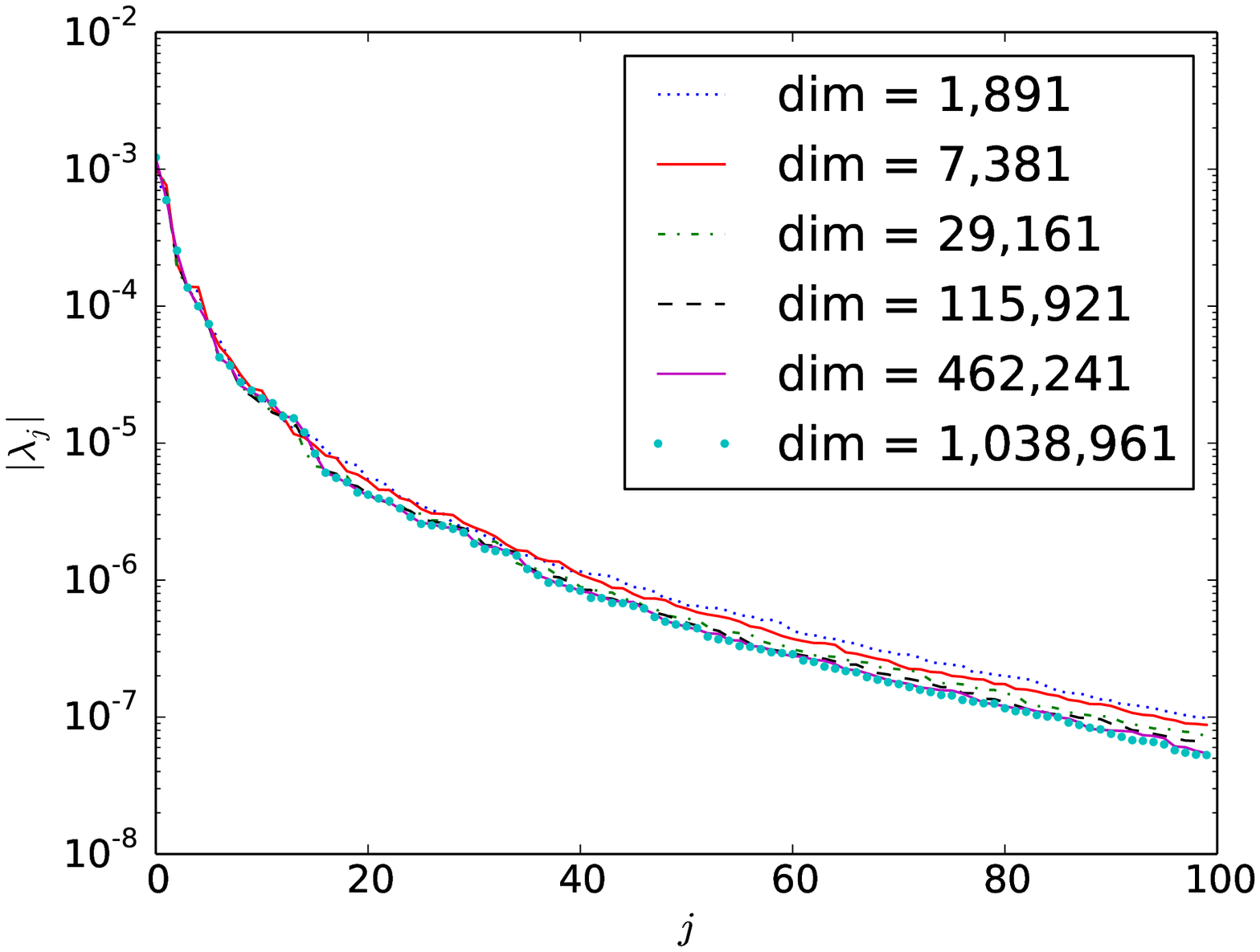}
\includegraphics[scale=0.31]{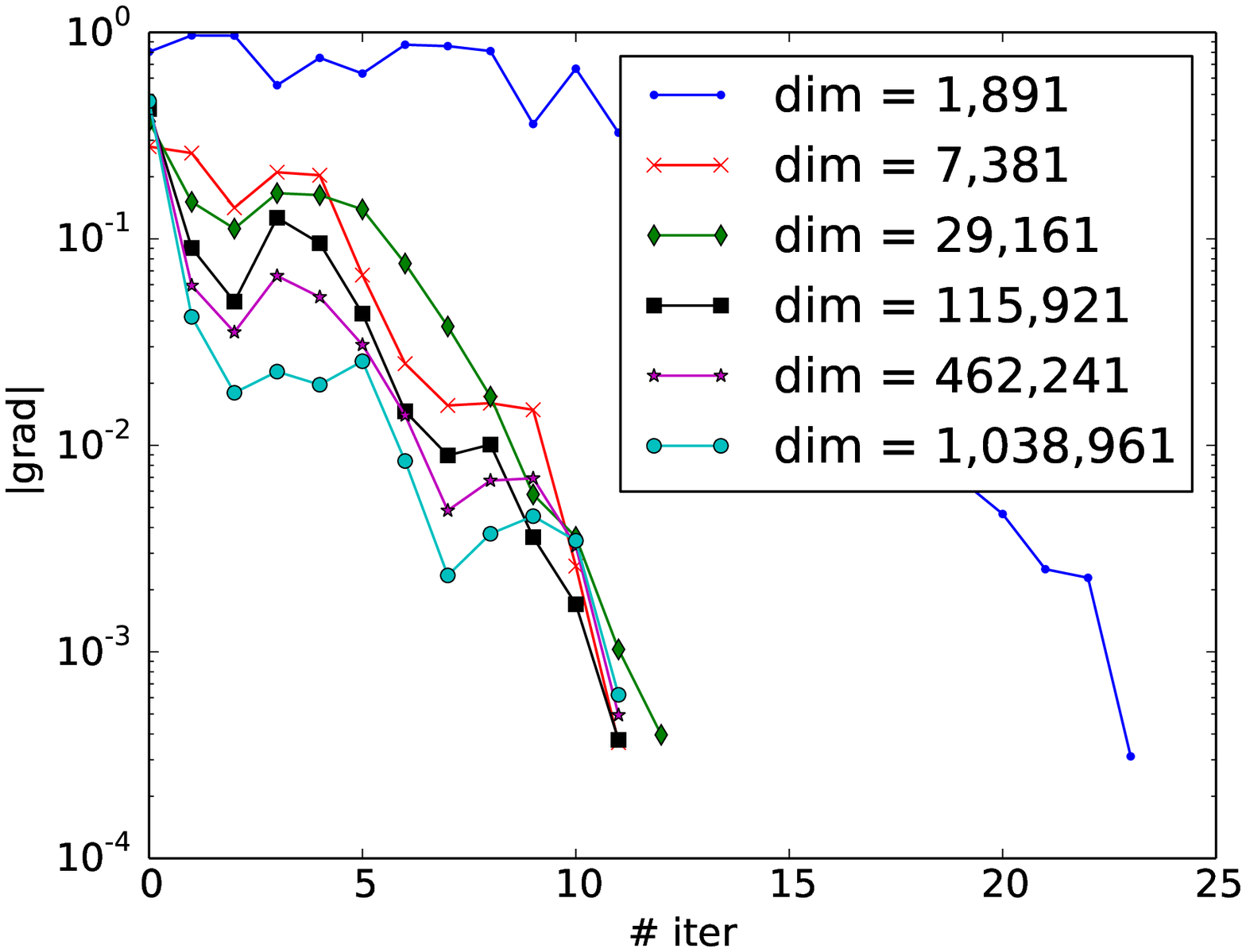}
\includegraphics[scale=0.31]{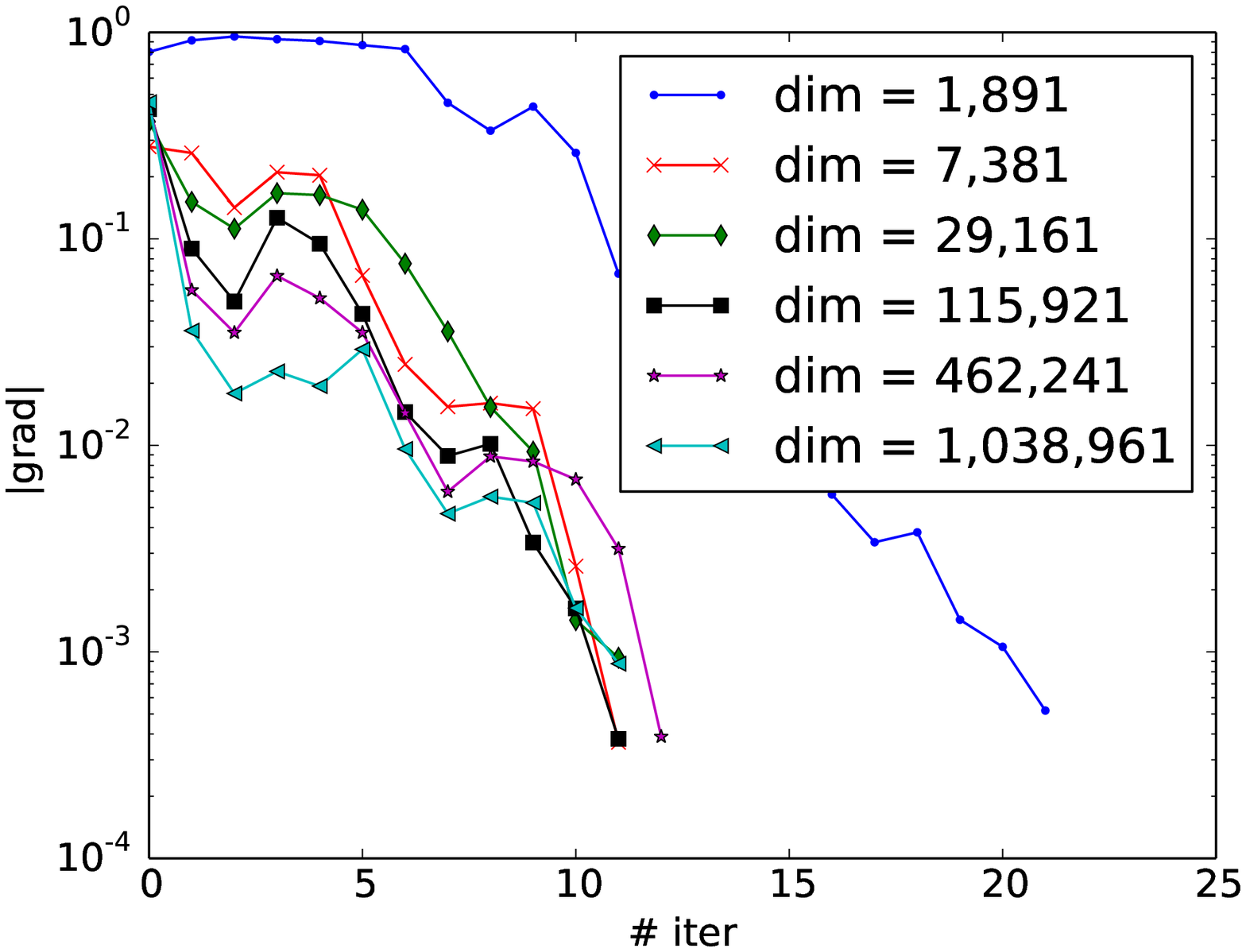}
\includegraphics[scale=0.31]{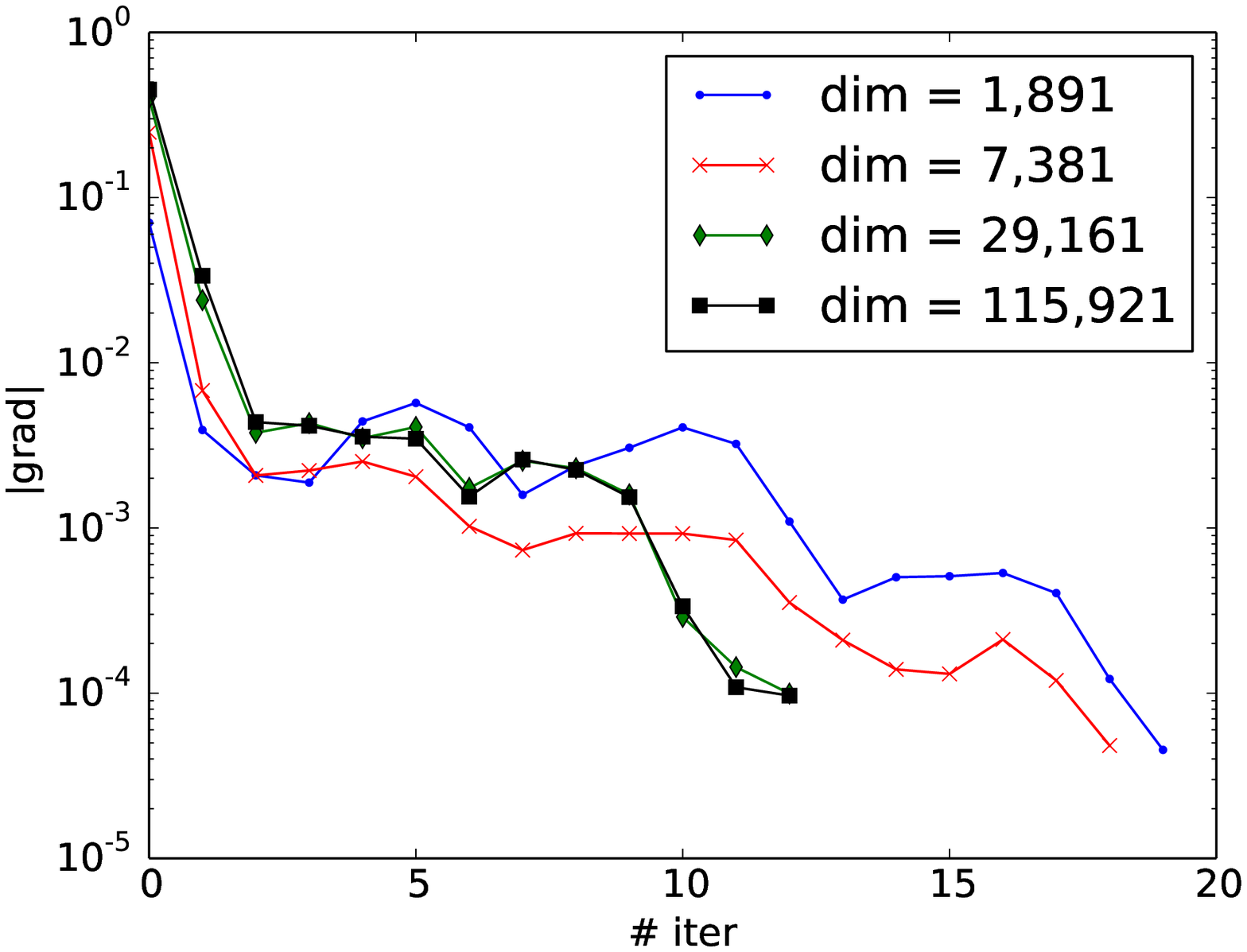}
\end{center}
\caption{Top-left: decay of the eigenvalues of the covariance-preconditioned Hessian (i.e., the generalized eigenvalues (in absolute value) of \eqref{eq:gEigen}) with different uncertain parameter dimensions (dim) at the optimal control obtained by linear approximation. Decay of the gradient with respect to the number of BFGS iterations (\# iter), with linear approximation (top-right), quadratic approximation (bottom-left), and linear approximation with Monte Carlo correction (bottom-right). Smaller dimensionality is shown for the last case due to constraints on the computating time.
}\label{fig:cost_mesh}
\end{figure}

The scalability of our approach with respect to the uncertain parameter dimension depends on three properties, namely the dimension-independent behavior of (\emph{i}) the spectral of the covariance preconditioned Hessian, (\emph{ii}) the convergence rate of BFGS, and (\emph{iii}) the variance of estimators for the Monte Carlo correction of the objective function. The top-left plot of Fig.\  \ref{fig:cost_mesh} demonstrates (\emph{i}): the eigenvalues of the covariance-preconditioned Hessian of the control objective (i.e., the generalized eigenvalues of problem \eqref{eq:gEigen}) at the optimal control obtained by linear approximation exhibit the same fast decay independent of the uncertain parameter dimension (ranging from a thousand to a million), thus indicating that approximating the trace by a randomized eigensolver is scalable with respect to the uncertain parameter dimension. The other three plots of Fig.\  \ref{fig:cost_mesh} provide numerical evidence of the mesh independent convergence of the BFGS algorithm, property (\emph{ii}), for this particular problem.  In the case of the linear and quadratic approximations using mesh continuation, the norm of the $z$-gradient drops below the prescribed tolerance within 12 iterations for all mesh resolutions except the coarsest one (which represents a pre-asymptotic result); in the case of the linear approximation with Monte Carlo correction,  BFGS converges within 19 iterations for all mesh resolutions using the solution of the linear approximation as initial guess. This implies that the gradient-based BFGS algorithm is also scalable with respect to the uncertain parameter dimension. Finally, Tables \ref{tab:VRmean_turbulence} and \ref{tab:VRvar_turbulence} indicate that the variance of the Monte Carlo correction , i.e., property (\emph{iii}), using Taylor approximation as control variates is independent of the uncertain parameter dimension, and, furthermore, it allows greater than 100X speed up with respect to plain Monte Carlo estimator to achieve the same accuracy (MSE). Combining the three dimension-independent properties, we conclude that the Taylor approximation and variance reduction stochastic optimization method is overall scalable with respect to the uncertain parameter dimension.

\section{Conclusions}
\label{sec:conclusion}

In this work we proposed a scalable computational framework to solve PDE-constrained optimal control under uncertainty where the cost functional involves the mean and variance of the control objective and the state equation is governed by a (nonlinear) PDE with high-dimensional uncertain parameter. Specifically, we advocate the use of (low-order) Taylor expansions of the control objective with respect to the uncertain parameter either to directly approximate the moments of the control objective or as control variates to reduce the variance of Monte Carlo estimator. The mean and variance expressions of the (linear and quadratic) Taylor expansions involve the trace of the covariance-preconditioned Hessian of the control objective with respect to the uncertain parameter.
Randomized algorithms for solution of generalized eigenvalue problems \cite{VillaPetraGhattas2017,SaibabaLeeKitanidis2016,SaibabaAlexanderianIpsen2016} allow for accurate and efficient approximation of this trace and only require computing the action of the covariance-preconditioned Hessian on a number of random directions that depends on its numerical rank. As we showed in the numerical results, this approach is more efficient and accurate than the Gaussian trace estimator when the eigenvalues of the covariance-preconditioned Hessian exhibit fast decay, which is true when the control objective is only sensitive to a limited number of directions in the uncertain parameter space. This is often true regardless of the complexity of the state problem.
Moreover, when the possibly biased Taylor (linear and quadratic) approximation is not sufficiently accurate, we can use it as a control variate and correct the mean and variance by Monte Carlo estimator of the remainder of the Taylor expansion, which was shown to achieve considerable computational cost reduction (by several orders of magnitude) compared to a plain Monte Carlo estimator. 
We have demonstrated the scalability of our method with respect to the uncertain parameter dimension in three aspects: the decay of the eigenvalues of the covariance-preconditioned Hessian of the control objective, the variance of the remainder of the Taylor expansion, and the number of optimization iterations. All three properties do not depend on the nominal dimension of the uncertain parameter but only on the intrinsic dimension of the uncertain parameter, i.e. on the number of directions the control objective is sensitive to. Moreover, by using the randomized eigensolver for trace estimator, Taylor approximation as control variate for variance reduction, and quasi-Newton method with a Lagrangian formalism for the optimization, we showed that the complexity---measured in number of PDE solves---scales linearly with respect to the intrinsic dimensionality to achieve target (high) accuracy, several orders of magnitude more accurate than a plain Monte Carlo estimator. This is illustrated numerically by solving an optimal boundary control governed by a nonlinear PDE model for turbulent jet flow with infinite-dimensional uncertain parameter field that characterizes the turbulent viscosity. Hence, the proposed computational framework for PDE-constrained optimal control under uncertainty is demonstrated scalable, accurate and efficient. 


Several extensions and further developments of the proposed computational framework are worth pursuing. First, we intend to investigate higher-order Taylor expansions beyond the quadratic approximation when the later becomes poor. Second, another research direction is to pursue more general risk measures than the mean-variance measure, such as VaR and conditional VaR \cite{KouriSurowiec2016, ShapiroDentchevaRuszczynski09}. Third, further computational savings can be achieved in combination with model reduction for the state PDE and the linearized PDEs \cite{ZahrFarhat2015, ChenQuarteroniRozza2016}. 
Fourth, one can replace the Monte Carlo estimator in the variance reduction by a sparse quadrature \cite{Chen2018} to possibly achieve faster convergence of the approximation error.
Fifth, application of the Taylor approximation and variance reduction in a multifidelity framework \cite{PeherstorferWillcoxGunzburger18} may lead to additional computational saving. Finally, one can apply the scalable computational framework to other challenging control and design problems, such as optimal design of metamaterials under uncertainty \cite{ChenGhattas18a, ChenGhattas18b}. 

\section*{Acknowledgement}
This work was supported by DARPA contract W911NF-15-2-0121, NSF grants CBET-1508713 and ACI-1550593, and DOE grants DE-SC0010518 and DE-SC0009286. We thank Robert Moser, Todd Oliver, and Myoungkyu Lee for the construction of the notional turbulent jet flow optimal boundary control under uncertainty problem.

\bibliographystyle{plain}
\bibliography{bibliography,ccgo}

\end{document}